\documentclass[12pt,a4paper,reqno]{article}
\usepackage[T1]{fontenc}
\usepackage{lmodern}
\usepackage[a4paper,left=25mm,top=28.4mm,textwidth=160mm,textheight=225mm,
headheight=0pt,headsep=0pt,footskip=12mm]{geometry}
\usepackage{amsmath,amssymb,amsthm,mathtools}
\usepackage{booktabs,array,graphicx,microtype,needspace}
\usepackage{indentfirst,xcolor,capt-of}
\usepackage[colorlinks=true,linkcolor=red,anchorcolor=green,
citecolor=blue,urlcolor=blue,bookmarksopen=true,
bookmarksnumbered=true,bookmarksdepth=2]{hyperref}
\hypersetup{pdftitle={Sharp Liouville theorems for the Hardy-H\'enon p-Laplace equation and classification of positive solutions for the critical Hardy p-Laplace equation},pdfauthor={Haoyang Lu; Zhitao Zhang},pdfsubject={Hardy--H\'enon p-Laplace equations, Liouville theorems and critical classification}}
\numberwithin{equation}{section}
\theoremstyle{plain}
\newtheorem{theorem}{Theorem}[section]
\newtheorem{corollary}[theorem]{Corollary}
\newtheorem{lemma}[theorem]{Lemma}

\newtheorem{remark}[theorem]{Remark}
\newcommand{\R}{\mathbb R}
\newcommand{\dd}{\,\mathrm d}
\newcommand{\norm}[1]{\lVert#1\rVert}

\allowdisplaybreaks[2]
\usepackage{enumitem,xurl}
\setlist[enumerate]{label=\textup{(\roman*)},leftmargin=*,itemsep=3pt}
\newcommand{\D}{D^{1,p}(\R^N)}
\newcommand{\eps}{\varepsilon}

\newcommand{\av}{\mathop{\mathrm{Av}}}
\DeclareMathOperator{\diver}{div}
\DeclareMathOperator{\supp}{supp}
\DeclareMathOperator{\Lip}{Lip}
\newcommand{\doi}[1]{\href{https://doi.org/#1}{doi:#1}}
\title{\large\textbf{Sharp Liouville theorems for the Hardy-H\'enon p-Laplace equation and classification of positive solutions for the critical Hardy p-Laplace equation}\thanks{
Z. Zhang is supported by the National Key R\&D Program of China
(2022YFA1005601) and the National Natural Science Foundation of China
(12671146,12031015).}}
\author{Haoyang Lu$^{a,b}$~~Zhitao Zhang$^{a,b}$\\
{\small $^{a}$ State Key Laboratory of Mathematical Sciences, Academy of Mathematics and Systems}\\
{\small Science, the Chinese Academy of Sciences, Beijing 100190, P. R. China}\\
{\small $^{b}$ School of Mathematical Sciences, University of Chinese Academy of Sciences,}\\
{\small Beijing 100049, P. R. China}\\
{\small E-mails: luhaoyang25@mails.ucas.ac.cn (H. Lu), zzt@math.ac.cn (Z. Zhang)}}
\date{}

\begin{document}
\maketitle
\begin{abstract}
We prove subcritical nonexistence for $-\Delta_pu=|x|^au^q$ in $\R^N$, with $1<p<N$ and $a>-p$, in the positive locally bounded weak class. Continuous Steiner symmetrization of a nonlinear auxiliary potential yields a directional pressure inequality; a localized Pohozaev identity and scaling give the Liouville theorem. At criticality, for $N\ge3$ and $-p<a<0$,
we classify all positive entire weak solutions in
$W^{1,p}_{\mathrm{loc}}(\R^N)$, without assuming local boundedness,
symmetry, stability, finite energy, or decay at infinity.
\\[0.3em]
\end{abstract}

\section{Introduction}\label{sec-intro}

We study nonexistence and critical classification for the equation
\begin{equation}\label{eq-scalar}
 -\Delta_pu=|x|^au^q\quad\text{in } \quad \R^N,
 \end{equation}
where $ \Delta_pu=\diver\bigl(|\nabla u|^{p-2}\nabla u\bigr), N\geq 2$, $1<p<N$, $a>-p$, and $q>p-1$. Scaling and supersolution estimates lead to the thresholds
\begin{align}
 q_H(a)&=\frac{p(N+a)}{N-p}-1,\label{eq-critical}\\
 q_S(a)&=\frac{(p-1)(N+a)}{N-p}.\label{eq-serrin-exponent}
\end{align}
Here
\begin{equation}\label{eq-threshold-gap}
 q_H(a)-q_S(a)=\frac{p+a}{N-p}>0.
\end{equation}
We prove nonexistence throughout $p-1<q<q_H(a)$, including the interval $q_S(a)<q<q_H(a)$ beyond the supersolution threshold. At $q=q_H(a)$ we also classify positive weak solutions for negative weights.

H\'enon's stellar model (1973) \cite{Henon1973} motivates $p=2$. Gidas and Spruck (1981) \cite{GS1981} proved, for nonnegative classical solutions,
\[
 N\ge3,\quad 1<q<\frac{N+2}{N-2},\quad
 -\Delta u=u^q\text{ in }\R^N
 \quad\Longrightarrow\quad u\equiv0.
\]
 Ni (1982) \cite{Ni1982} initiated weighted Dirichlet existence theory on the unit ball.

For the quasilinear equation, Bidaut-V\'eron (1989) \cite{BV1989} studied radial solutions and singularities. Kawano, Ni, and Yotsutani (1990) \cite{KNY1990} developed the radial Pohozaev method. The radial theory gives, for $a\ge0$, the criterion recalled in \cite[Proposition~A]{PD}:
\[
 \text{a positive regular radial entire solution exists}
 \quad\Longleftrightarrow\quad q\ge q_H(a).
\]

Bidaut-V\'eron and Garc\'ia-Huidobro (2001) \cite{BGH2001} treated weighted singularities and nonradial estimates. For $a\ge0$ and $1<p<N$, the integral methods of
Bidaut-V\'eron--Pohozaev \cite{BidautPohozaev2001}
and Mitidieri--Pohozaev \cite{MP2001} yield
\[
 p-1<q\le\frac{(p-1)(N+a)}{N-p}
 \quad\Longrightarrow\quad
 \nexists\,u>0:\;
 -\Delta_pu\ge |x|^a u^q
 \quad\text{in }\mathbb R^N.
\]
See also \cite[Eq.~(1.4)]{PD}.

Serrin and Zou (2002) \cite{SZ} established unrestricted nonexistence and universal local bounds for
\[
 a=0,\qquad p-1<q<q_H(0)=\frac{Np}{N-p}-1.
\]
Pol\'a\v cik, Quittner, and Souplet (2007) \cite{PQS2007} developed doubling and rescaling to derive singularity and decay estimates from such Liouville theorems.

Pucci and Servadei (2008) \cite{PS} studied radial ground states with singular weights and supplied the critical profile used below. Their 2010 Pohozaev results \cite{PS2010} concern weak solutions in global energy spaces.

Phan and Souplet (2012) \cite{PhanSouplet2012} proved full positive-weight subcritical nonexistence for bounded solutions when $p=2$, $N=3$, and also treated negative weights and decay estimates.

For $a\ge0$, Phan and Duong (2015) \cite[Theorems~1.2--1.3]{PD} proved
\begin{align}
 N<p+1,\quad p-1<q<q_H(a)
 &\quad\Longrightarrow\quad \text{nonexistence},\label{eq-PD-lowdim}\\
 p-1<q<q_H(0)
 &\quad\Longrightarrow\quad \text{nonexistence in every }N>p.\label{eq-PD-unweighted-range}
\end{align}
No global bound is assumed \cite[Remark~1.4]{PD}. The proofs use spherical Sobolev feedback and decay, respectively. Their Conjecture~B asks for the full range $q<q_H(a)$. For $a\ge0$, regularity identifies their $C^1$ weak class with our locally bounded weak class.

For $p=2$, Garc\'ia-Meli\'an's 2017 preprint \cite{GarciaMelian2017} presents full subcritical nonexistence for $a>-2$ in $H^1_{\mathrm{loc}}\cap L^\infty_{\mathrm{loc}}$, using Emden--Fowler monotonicity and derived stability. Li and Zhang (2019) \cite{LiZhang2019Monotonicity} established a cylinder monotonicity theorem and derived sharp existence criteria for weighted semilinear equations with positive solutions bounded near the origin.

Ciraolo, Figalli, and Roncoroni (2020) \cite{CFR2020} used differential identities to classify finite-energy critical solutions for anisotropic $p$-Laplace operators in convex cones. For the semilinear equation, Giga and Ng\^o (2022) \cite{GigaNgo2022} classify real powers and weights, imposing continuity at zero for singular weights.

At criticality, Ciraolo and Corso (2022) \cite{CiraoloCorso2022} used a weighted Riemannian formulation to prove symmetry results for Caffarelli--Kohn--Nirenberg equations. Le and Le (2023) \cite{LeLe2023} classified positive $D^{1,p}(\R^N)$ solutions for $-p<a<0$. For the unweighted critical equation, Zhang
\cite[Corollary~1.4]{Zhang2026Critical}
classified all positive entire weak solutions in
$W^{1,p}_{\mathrm{loc}}(\R^N)$ for $1<p<N$,
without additional global assumptions.
Le \cite[Theorem~1.1]{Le2026Anisotropic}
obtained the corresponding classification for anisotropic
$p$-Laplace operators.
Theorem~\ref{thm-negative-critical} treats the singular-weight range
$-p<a<0$ for $N\ge3$ in the same local weak class,
without assuming local boundedness, symmetry, stability,
finite energy, or decay at infinity.
Local boundedness, radial symmetry, and finite energy
are conclusions of the theorem.

Du and Zhou (2024) \cite[Sections~3--4]{DZ} proved radial subcritical nonexistence and nonradial nonexistence under
\begin{equation}\label{eq-prior-decay}
 u(x)\le C|x|^{-(p+a)/(q-p+1)},\qquad x\ne0.
\end{equation}
For negative weights at criticality, Pan and Shang (2026) \cite{PanShang2026} obtained classifications under energy-growth or asymptotic assumptions.  For positive weights, Dai, Duan, Gui, and Li (2026) \cite{DDGL2026} constructed nonradial solutions at specific bifurcation parameters. Thus a classification by radial profiles cannot hold for all positive weights.

For the subcritical theorem, the solution class is
\begin{equation}\label{eq-solution-class}
 u\in W^{1,p}_{\mathrm{loc}}(\R^N)\cap L^\infty_{\mathrm{loc}}(\R^N),
 \qquad u\ge0\quad\text{a.e.}
\end{equation}
with
\begin{equation}\label{eq-weak}
 \int_{\R^N}|\nabla u|^{p-2}\nabla u\cdot\nabla\phi\dd x
 =\int_{\R^N}|x|^au^q\phi\dd x
 \quad\text{for all }\phi\in C_c^\infty(\R^N).
\end{equation}
Here $u\in L^\infty_{\mathrm{loc}}$ and $a>-p>-N$ imply $|x|^au^q\in L^1_{\mathrm{loc}}$; \eqref{eq-weak} includes the origin. We also use
\[
 \D=\overline{C_c^\infty(\R^N)}^{\,\norm{\nabla\cdot}_{L^p}},
 \qquad \norm\phi_{\D}=\norm{\nabla\phi}_{L^p}.
\]

\begin{theorem}[Subcritical Liouville theorem]\label{thm-scalar}
Let $1<p<N$, $a>-p$, and $p-1<q<q_H(a)$. Every nonnegative weak entire solution of \eqref{eq-scalar} satisfying
\[
 u\in W^{1,p}_{\mathrm{loc}}(\R^N)\cap L^\infty_{\mathrm{loc}}(\R^N)
\]
is zero almost everywhere in $\R^N$.
\end{theorem}

\Needspace{9\baselineskip}
\begin{corollary}[$C^{1,\alpha}$ solutions away from the origin]\label{cor-classical}
Under the parameter assumptions of Theorem~\ref{thm-scalar}, let $u\ge0$ satisfy, for some $\alpha\in(0,1)$,
\begin{equation}\label{eq-classical-class}
 u\in C^{1,\alpha}_{\mathrm{loc}}(\R^N\setminus\{0\})\cap C(\R^N)
\end{equation}
and $-\Delta_pu=|x|^au^q$ weakly on $\R^N\setminus\{0\}$. Then $u\equiv0$. In particular, there is no positive solution in this class.
\end{corollary}
In Corollary \ref{cor-classical} the equation is understood in the weak sense:
\begin{equation}\label{eq-punctured-weak}
 \int_{\R^N}|\nabla u|^{p-2}\nabla u\cdot\nabla\phi\dd x
 =\int_{\R^N}|x|^au^q\phi\dd x,
 \quad \forall\phi\in C_c^\infty(\R^N\setminus\{0\}).
\end{equation}
This includes zero-gradient points when $1<p<2$. The origin is removable by Section~\ref{sec-origin}.

 For $a\ge0$, it gives the conclusion of \cite[Conjecture~B]{PD} for every $N>p$. In particular, it includes the parameters outside \eqref{eq-PD-lowdim}, \eqref{eq-PD-unweighted-range}, and the supersolution range:
\begin{equation}\label{eq-new-comparison-range}
\begin{gathered}
 a>0,\qquad N\ge p+1,\qquad p\ne2,\\
 q\ge q_H(0),\qquad q>q_S(a),\qquad q<q_H(a).
\end{gathered}
\end{equation}

The profile \eqref{eq-bubble} proves sharpness at $q=q_H(a)$ in both classes.\par
For classification of positive solutions of \eqref{eq-scalar} for critical exponent with negative weights, we have:
\begin{theorem}\label{thm-negative-critical}
Let $N\ge3$, $1<p<N$, and $-p<a<0$. Every positive weak entire solution
\[
 u\in W^{1,p}_{\mathrm{loc}}(\R^N),\qquad
 -\Delta_pu=|x|^au^{q_H(a)}\quad\text{in }\R^N,
\]
where the weak identity includes the origin, has the form
\begin{equation}\label{eq-negative-classification}
 u(x)=\lambda^{(N-p)/p}U(\lambda x),\qquad \lambda>0,
\end{equation}
where
\[
\begin{gathered}
 U(x)=C_0\left(1+|x|^{\frac{p+a}{p-1}}\right)^{-\frac{N-p}{p+a}},\\
 C_0=\left[(N+a)\left(\frac{N-p}{p-1}\right)^{p-1}\right]^{\frac{N-p}{p(p+a)}}.
\end{gathered}
\]
In particular,
\begin{equation}\label{eq-negative-energy}
 u\in\D,\qquad
 \int_{\R^N}|\nabla u|^p\dd x
 =\int_{\R^N}|x|^au^{q_H(a)+1}\dd x<\infty.
\end{equation}
\end{theorem}

In Theorem~\ref{thm-negative-critical}, local boundedness, radial symmetry, and finite energy are conclusions. Positive solutions in \eqref{eq-classical-class} belong to this weak class after removal of the origin. The proof in Section~\ref{sec-negative-critical} first obtains a pressure bound on a weighted cone, then compares $u$ with an explicit critical profile by Picone's inequality.

For the variational interpretation, set
\[
 m_a=q_H(a)+1=\frac{p(N+a)}{N-p},\qquad -p<a<0.
\]
The weighted Lebesgue space and its norm are
\[
\begin{gathered}
 L^{m_a}(\R^N,|x|^a\dd x)
 =\{w:\R^N\to\R\text{ measurable}:\norm{w}_{m_a,a}<\infty\},\\
 \norm{w}_{m_a,a}
 =\left(\int_{\R^N}|x|^a|w|^{m_a}\dd x\right)^{1/m_a}.
\end{gathered}
\]
The Hardy--Sobolev embedding is
\begin{equation}\label{eq-critical-hardy-sobolev}
\begin{gathered}
 \D\hookrightarrow L^{m_a}(\R^N,|x|^a\dd x),\\
 S_{p,a}\norm{w}_{m_a,a}^{p}
 \le\int_{\R^N}|\nabla w|^p\dd x,\qquad w\in\D,
\end{gathered}
\end{equation}
where the optimal constant is
\[
 S_{p,a}:=\inf_{0\ne w\in\D}
 \frac{\displaystyle\int_{\R^N}|\nabla w|^p\dd x}
 {\left(\displaystyle\int_{\R^N}|x|^a|w|^{m_a}\dd x\right)^{p/m_a}}>0.
\]
The source space $\D$ is defined by the unweighted gradient norm above; the weight enters the target measure. The sharp inequality and its extremal functions are given in \cite{GhoussoubYuan2000}.

\Needspace{12\baselineskip}
\begin{corollary}[Hardy--Sobolev minimizers and common energy]\label{cor-critical-minimizers}
Let $N\ge3$, $1<p<N$, and $-p<a<0$. Every positive weak entire solution
\[
 u\in W^{1,p}_{\mathrm{loc}}(\R^N),\qquad
 -\Delta_pu=|x|^au^{m_a-1}\quad\text{in }\R^N,
\]
with the weak identity including the origin, belongs to $\D$ and attains $S_{p,a}$:
\begin{equation}\label{eq-critical-minimality}
 \frac{\displaystyle\int_{\R^N}|\nabla u|^p\dd x}
 {\left(\displaystyle\int_{\R^N}|x|^au^{m_a}\dd x\right)^{p/m_a}}
 =S_{p,a}.
\end{equation}
Moreover, all these solutions have the same energy:
\begin{equation}\label{eq-critical-common-energy}
 \int_{\R^N}|\nabla u|^p\dd x
 =\int_{\R^N}|x|^au^{m_a}\dd x
 =S_{p,a}^{\frac{N+a}{p+a}}.
\end{equation}
\end{corollary}
\begin{proof}
The profile $U$ in Theorem~\ref{thm-negative-critical} is an extremal for \eqref{eq-critical-hardy-sobolev} \cite{GhoussoubYuan2000}. Critical dilations preserve the Sobolev quotient, so Theorem~\ref{thm-negative-critical} gives $u\in\D$ and \eqref{eq-critical-minimality}. By \eqref{eq-negative-energy}, put
\[
 T:=\int_{\R^N}|\nabla u|^p\dd x
 =\int_{\R^N}|x|^au^{m_a}\dd x.
\]
Then
\[
 S_{p,a}=T^{1-p/m_a},\qquad
 T=S_{p,a}^{m_a/(m_a-p)}
   =S_{p,a}^{(N+a)/(p+a)},
\]
which proves \eqref{eq-critical-common-energy}.
\end{proof}

\medskip\noindent\textbf{The obstruction to linear-potential arguments.}
OpenAI's virial and source-layer method \cite[Sections~2--5]{OpenAI2026} was adapted to fractional kernels by Mei and Wei \cite{MW2026}. Li, Lu, and Zhang \cite{LiLuZhang2026Fractional} developed localized kernel estimates for weighted fractional H\'enon--Lane--Emden systems, proving subcritical nonexistence and examples of both existence and nonexistence on one critical hyperbola. Deng and Zhang \cite{DengZhang2026} treated variable kernels and coefficients. These arguments use a linear potential representation. The obstruction for $p\ne2$ is already visible in the following pairing identity.

For $N\ge3$, let $\Phi(x)=|x|^{2-N}/((N-2)|\mathbb S^{N-1}|)$ and $T_2g=\Phi*g$. For $0\le f_1,f_2\in C_c^\infty(\R^N)$,
\begin{equation}\label{eq-intro-linear-reciprocity}
 \int f_1T_2f_2
 =\iint f_1(x)f_2(y)\Phi(x-y)\dd x\dd y
 =\int f_2T_2f_1.
\end{equation}
The shared kernel gives the matching interaction terms and pressure cancellation in \cite[Eq.~(3.1)]{OpenAI2026}.

For $0\le g\in C_c^\infty$, let $T_pg\in\D$ be the unique weak solution of $-\Delta_p(T_pg)=g$. Then
\begin{equation}\label{eq-intro-nonlinear-homogeneity}
 T_p(cg)=c^{1/(p-1)}T_pg,\qquad c>0.
\end{equation}
For $p\ne2$, $g\not\equiv0$, and $c>0$, $c\ne1$, put $E_g=\int gT_pg>0$. Then
\begin{equation}\label{eq-intro-failed-reciprocity}
 \int(cg)T_pg=cE_g,
 \qquad
 \int gT_p(cg)=c^{1/(p-1)}E_g\ne cE_g.
\end{equation}
Thus linear reciprocity fails even for proportional sources. Two-sided Wolff bounds do not recover it.

\medskip\noindent\textbf{The auxiliary potential and directional pressure.}
Lu and Zhang \cite{LuZhang2026Disk} combined divergence inequalities, interpolation, and scaling to obtain a disk-shaped nonexistence region for the Lane--Emden system. For the present scalar equation, the estimate needed for scaling is a local bound uniform in the solution. We obtain it by varying the source of a nonlinear auxiliary potential.

For $f=|x|^au^q$ and an annular cutoff $\psi$, set
\[
 g=\psi^{p-1}f,\qquad z=T_pg,\qquad v=\psi u.
\]
Here $z\in\D$ although $u$ need not have finite energy. For a source $g$, define the dual functional
\begin{equation}\label{eq-intro-dual-functional}
 W_p(g)=\sup_{\phi\in\D}
 \left\{\int g\phi-\frac1p\int|\nabla\phi|^p\right\}.
\end{equation}
Let $(g_\tau)_{\tau\ge0}$ be the continuous Steiner symmetrization of $g$ in a fixed direction. Brock's inequalities \cite{Brock1995,Brock2000,Brock} give $W_p(g_\tau)\ge W_p(g)$. Differentiating at $\tau=0$ yields the directional pressure inequality in Lemma~\ref{lem-pressure}.

Picone and convexity bound $z-\psi u$ in gradient and weighted $L^p$ norms. Wolff bounds then control the source-layer lengths that enter the endpoint errors. The inequalities
\[
 \frac1p-\frac1{q+1}>0,\qquad
 \frac{N+a}{q+1}-\frac{N-p}{p}>0
\]
give, through the localized Pohozaev identity and iteration,
\[
 \int_{B_1}|x|^au^{q+1}\dd x\le C(N,p,a,q)
\]
uniformly in $u$. The error bound has a fixed sublinear power of the local reaction energy, which permits absorption. Scaling then proves Theorem~\ref{thm-scalar}.

The critical proof uses a different pressure function. For $a<0$, the critical power is strictly subcritical for the unweighted equation, giving coarse decay by rescaling. In a weighted-cone metric, a differential identity and logarithmic cutoffs yield a pressure upper bound. A local flux comparison locates the maximum of $u$ at the origin, and Picone comparison with an explicit profile gives both classification and finite energy.

Sections~\ref{sec-basic}--\ref{sec-close} develop the auxiliary-potential estimates, the directional pressure inequality, and the subcritical proof. Section~\ref{sec-negative-critical} contains the critical classification.

\section{Source mass and the nonlinear auxiliary potential}\label{sec-basic}

In Sections~\ref{sec-basic}--\ref{sec-close}, assume $1<p<N$, $a>-p$, $p-1<q<q_H(a)$, and \eqref{eq-solution-class}--\eqref{eq-weak}. We first justify the weak tests and removal of the origin, then estimate the source mass and compare an auxiliary potential with the cutoff solution. The different weak class in Theorem~\ref{thm-negative-critical} is treated separately in Section~\ref{sec-negative-critical}.

\subsection{Regularity of weak solutions}\label{sec-regularity}

Put $f=|x|^au^q$. For some $\alpha\in(0,1)$, \cite[Theorem~2]{DiBenedetto1983} and \cite[Theorem~1]{Tolksdorf} give
\[
 a\ge0\quad\Longrightarrow\quad f\in L^\infty_{\mathrm{loc}}(\R^N)
 \quad\Longrightarrow\quad u\in C^{1,\alpha}_{\mathrm{loc}}(\R^N).
\]
For $-p<a<0$, \cite[Theorem~8]{Serrin64} and \cite[Theorem~1]{Tolksdorf} give
\[
\begin{gathered}
 \frac Np<r<\frac N{-a},\qquad f\in L^r_{\mathrm{loc}}(\R^N)
 \quad\Longrightarrow\quad u\in C(\R^N),\\
 f\in L^\infty_{\mathrm{loc}}(\R^N\setminus\{0\})
 \quad\Longrightarrow\quad u\in C^{1,\alpha}_{\mathrm{loc}}(\R^N\setminus\{0\}).
\end{gathered}
\]
For this representative, the strong minimum principle follows from Serrin's Harnack inequality \cite[Theorem~5]{Serrin64} (see also \cite[Theorem~1.14.8]{zzt}):
\[
 -\Delta_pu=f\ge0,\quad u\ge0
 \quad\Longrightarrow\quad u\equiv0\ \text{or}\ u>0\text{ in }\R^N.
\]

For $u>0$, the source satisfies $f\in C^1(\R^N\setminus\{0\})$. To obtain Pohozaev's identity for weak solutions, we use the Pucci--Serrin domain-variation identity \cite[Lemma~1]{DMS}: the required local Lipschitz regularity follows from the above $C^{1,\alpha}$ estimates away from $0$. We then choose a radial vector field and remove the cutoff at the origin.

\subsection{Removal of the origin}\label{sec-origin}

We prove $u\in W^{1,p}_{\mathrm{loc}}(\R^N)$ and \eqref{eq-weak} for Corollary~\ref{cor-classical}.

\begin{proof}[Proof of Corollary~\ref{cor-classical}]
For $u$ as in the corollary and $f=|x|^au^q$,
\[
 a>-N,\quad u\in C(\R^N)\quad\Longrightarrow\quad
 u\in L^\infty_{\mathrm{loc}},\quad f,fu\in L^1_{\mathrm{loc}}.
\]

Fix $R>0$, $M_R=\max_{\overline{B_{2R}}}u$, and smooth $0\le\zeta,\chi_\eps\le1$ satisfying
\[
\begin{gathered}
 \zeta\in C_c^\infty(B_{2R}),\qquad
 \zeta=1\text{ on }B_R,\qquad |\nabla\zeta|\le C/R,\\
 \chi_\eps=0\text{ on }B_\eps,\qquad
 \chi_\eps=1\text{ outside }B_{2\eps},\qquad
 |\nabla\chi_\eps|\le C/\eps.
\end{gathered}
\]
Put $\eta_\eps=\zeta\chi_\eps$, $0<\eps<R/4$. Since $u\in C^1_{\mathrm{loc}}(\R^N\setminus\{0\})$, the test $u\eta_\eps^p\in C_c^1(\R^N\setminus\{0\})$ is admissible in \eqref{eq-punctured-weak} by smooth approximation. Testing and Young's inequality give
\[
\begin{aligned}
 \int\eta_\eps^p|\nabla u|^p
 &=\int fu\eta_\eps^p
   -p\int u\eta_\eps^{p-1}|\nabla u|^{p-2}\nabla u\cdot\nabla\eta_\eps\\
 &\le\int fu\eta_\eps^p
   +\frac12\int\eta_\eps^p|\nabla u|^p
   +C\int u^p|\nabla\eta_\eps|^p.
\end{aligned}
\]
Hence, with $C$ independent of $\eps$,
\begin{equation}\label{eq-origin-energy}
\begin{aligned}
 \int_{B_R\setminus B_{2\eps}}|\nabla u|^p
 &\le C\int u^p|\nabla(\zeta\chi_\eps)|^p
       +C\int fu(\zeta\chi_\eps)^p\\
 &\le CM_R^p(R^{N-p}+\eps^{N-p})
       +CM_R^{q+1}R^{N+a}.
\end{aligned}
\end{equation}
Fatou's lemma in \eqref{eq-origin-energy} gives $\nabla u\in L^p(B_R)$. Also,
\[
 \norm{u\nabla\chi_\eps}_{L^p(B_R)}^p
 \le CM_R^p\eps^{N-p}\longrightarrow0.
\]
Hence
\[
 \chi_\eps u\longrightarrow u\text{ in }L^p(B_R),\qquad
 \nabla(\chi_\eps u)=\chi_\eps\nabla u+u\nabla\chi_\eps
 \longrightarrow\nabla u\text{ in }L^p(B_R).
\]
The closedness of weak differentiation identifies this $L^p$ limit as the weak gradient of $u$. Thus $u\in W^{1,p}_{\mathrm{loc}}(\R^N)$.

Test \eqref{eq-punctured-weak} with $\phi\chi_\eps$, $\phi\in C_c^\infty(\R^N)$:
\begin{equation}\label{eq-origin-flux}
\begin{aligned}
 \left|\int\phi|\nabla u|^{p-2}\nabla u\cdot\nabla\chi_\eps\right|
 &\le C\norm\phi_\infty
 \left(\int_{B_{2\eps}\setminus B_\eps}|\nabla u|^p\right)^{(p-1)/p}
 \eps^{(N-p)/p}\\
 &\longrightarrow0.
\end{aligned}
\end{equation}
The other two terms converge by dominated convergence, since $f\in L^1_{\mathrm{loc}}$ and $|\nabla u|^{p-1}\in L^{p/(p-1)}_{\mathrm{loc}}$. Thus \eqref{eq-weak} holds. Theorem~\ref{thm-scalar} and continuity yield $u\equiv0$. The removal argument itself uses only $1<p<N$, $a>-N$, and $q>0$; it will also apply to the critical equation.
\end{proof}

\subsection{Source mass and Picone inequality}

Henceforth, assume $u>0$, by Section~\ref{sec-regularity}.

Set
\begin{equation}\label{eq-basic-notation}
\begin{gathered}
 m=q+1>p,\qquad d=q-p+1=m-p>0,\qquad d_0=\frac1p-\frac1m>0,\\
 f=|x|^au^q,\qquad V=\frac{f}{u^{p-1}},\qquad
 F(R)=\int_{B_R}fu\dd x,\\
 \Gamma=\frac{N+a}{m}-\frac{N-p}{p}>0.
\end{gathered}
\end{equation}
Here $F(R)<\infty$ for $R<\infty$, and $\D\hookrightarrow L^{p^*}$, $p^*=Np/(N-p)$. Unspecified spatial integrals are over $\R^N$.

The source-mass bound below controls the lengths of superlevel intervals in Section~\ref{sec-pressure}.

\begin{lemma}\label{lem-basic}
There is a constant $C=C(N,p)$ such that, for every $x\in\R^N$ and $\ell>0$,
\begin{equation}\label{eq-source-ball}
 \int_{B_\ell(x)} f\dd y\le C\ell^{N-p}u(x)^{p-1}.
\end{equation}
Moreover, every $\phi\in\D$ satisfies
\begin{equation}\label{eq-picone}
 \int_{\R^N}V|\phi|^p\dd x\le\int_{\R^N}|\nabla\phi|^p\dd x.
\end{equation}
\end{lemma}
\begin{proof}
For $\mu=f\dd x$, the local Wolff bound \cite[Theorem~1.6]{KM1994}, also stated in \cite[Theorem~2.6]{Kilpelainen}, gives
\[
 u(x)\ge c\int_\ell^{2\ell}
 \left(\frac{\mu(B_t(x))}{t^{N-p}}\right)^{1/(p-1)}\frac{\dd t}{t}
 \ge c\left(\frac{\mu(B_\ell(x))}{\ell^{N-p}}\right)^{1/(p-1)}.
\]

For $0\not\equiv\phi\in C_c^\infty(\R^N)$, set
\[
 w=|\phi|,\qquad
 \min_{\supp\phi}u>0,\qquad
 \psi=\frac{w^p}{u^{p-1}}\in W^{1,p}(\R^N)\cap L^\infty(\R^N),
 \qquad \supp\psi\subset\supp\phi.
\]
Bounded mollification makes $\psi$ admissible in \eqref{eq-weak}, since $f\in L^1_{\mathrm{loc}}$. Picone's inequality \cite{AH1998} gives
\[
 |\nabla w|^p-
 |\nabla u|^{p-2}\nabla u\cdot\nabla\psi\ge0.
\]
Since $|\nabla w|\le|\nabla\phi|$,
\[
 \int V|\phi|^p=\int f\psi
 =\int |\nabla u|^{p-2}\nabla u\cdot\nabla\psi
 \le\int|\nabla\phi|^p.
\]
For $\phi\in\D$, density and Sobolev give $\phi_j\in C_c^\infty$ such that
\[
 \norm{\nabla(\phi_j-\phi)}_{L^p}\longrightarrow0,\qquad
 \norm{\phi_j-\phi}_{L^{p^*}}\longrightarrow0,\qquad
 \phi_j\longrightarrow\phi\quad\text{a.e. along a subsequence}.
\]
Fatou's lemma and $V\ge0$ yield
\[
 \int V|\phi|^p
 \le\liminf_{j\to\infty}\int V|\phi_j|^p
 \le\lim_{j\to\infty}\int|\nabla\phi_j|^p
 =\int|\nabla\phi|^p.
\]
This proof also applies at criticality after local regularity is established.
\end{proof}

\subsection{Potential comparison and convexity remainders}

We compare the cutoff solution $\psi u$ with the potential of $\psi^{p-1}f$.

Let $(p^*)'=p^*/(p^*-1)$. For $0\le g\in L^{(p^*)'}(\R^N)$, define $T_pg=\operatorname*{argmin}_{\phi\in\D}J_g(\phi)$, where
\begin{equation}\label{eq-auxiliary-functional}
 J_g(\phi)=\frac1p\int|\nabla\phi|^p\dd x-\int g\phi\dd x.
\end{equation}
Sobolev gives
\[
 \left|\int g\phi\right|
 \le\norm g_{L^{(p^*)'}}\norm\phi_{L^{p^*}}
 \le C\norm g_{L^{(p^*)'}}\norm{\nabla\phi}_{L^p}.
\]
In particular,
\[
 J_g(\phi)\ge\frac1p\norm\phi_{\D}^p
 -C\norm g_{L^{(p^*)'}}\norm\phi_{\D}
 \longrightarrow\infty\quad\text{as }\norm\phi_{\D}\to\infty.
\]
The direct method and strict convexity give a unique minimizer $z$. Testing its Euler equation with $z_-=\max\{-z,0\}$ and $z$ gives
\[
\begin{gathered}
 -\int|\nabla z_-|^p=\int gz_-\ge0
 \quad\Longrightarrow\quad z\ge0,\\
 -\Delta_pz=g,\qquad \int|\nabla z|^p=\int gz<\infty.
\end{gathered}
\]

Now choose an annular cutoff $\psi\ge0$ whose positive powers are $C^1$. Set
\begin{equation}\label{eq-potential-notation}
\begin{gathered}
 g=\psi^{p-1}f\in C_c^1(\R^N),\qquad
 v=\psi u,\qquad z=T_pg,\qquad h=z-v,\\
 F_\psi=\int\psi^pfu\dd x=\int Vv^p\dd x=\int gv\dd x.
\end{gathered}
\end{equation}
In particular, $g=Vv^{p-1}$.

\Needspace{11\baselineskip}
\begin{lemma}\label{lem-potential}
Put $K=\norm\psi_\infty$, $E_z=\int|\nabla z|^p$, and $E_v=\int|\nabla v|^p$. Then
\begin{equation}\label{eq-pointwise-potential}
 0\le z\le Ku,\qquad |h|\le Ku,\qquad E_z\le F_\psi.
\end{equation}
Moreover, $R_\psi:=E_v-F_\psi\ge0$, and
\begin{equation}\label{eq-h-energy}
 \int|\nabla h|^p\dd x\le
 \begin{cases}
 CR_\psi,&p\ge2,\\
 CR_\psi^{p/2}(F_\psi+R_\psi)^{(2-p)/2},&1<p<2.
 \end{cases}
\end{equation}
The same right-hand side bounds $\int V|h|^p$.
\end{lemma}
\begin{proof}

Let $z_R\in W_0^{1,p}(B_R)$ solve $-\Delta_pz_R=g$ and extend it by zero. Comparison gives
\[
 g\le K^{p-1}f=-\Delta_p(Ku),\quad z_R|_{\partial B_R}=0\le Ku
 \quad\Longrightarrow\quad 0\le z_R\le Ku.
\]
Density of $C_c^\infty$ in $\D$ gives
\[
 J_g(z_R)=\min_{W_0^{1,p}(B_R)}J_g\downarrow\min_{\D}J_g=J_g(z),
 \qquad
 \int|\nabla z_R|^p=-\frac p{p-1}J_g(z_R)\longrightarrow\int|\nabla z|^p.
\]
Weak compactness, uniqueness, and uniform convexity give
\[
 z_R\rightharpoonup z\text{ in }\D,\qquad
 \norm{z_R}_{\D}\longrightarrow\norm z_{\D}
 \quad\Longrightarrow\quad z_R\longrightarrow z\text{ in }\D.
\]
Passing to an a.e. convergent subsequence,
\[
 0\le z,v\le Ku\quad\Longrightarrow\quad |h|=|z-v|\le Ku.
\]
Hence we prove the pointwise comparison.

By H\"older and \eqref{eq-picone},
\begin{equation}\label{eq-potential-energy-chain}
\begin{aligned}
 E_z=\int Vv^{p-1}z
 &\le\left(\int Vv^p\right)^{(p-1)/p}
       \left(\int Vz^p\right)^{1/p}\\
 &\le F_\psi^{(p-1)/p}E_z^{1/p}.
\end{aligned}
\end{equation}
Thus $E_z\le F_\psi$.

Finally, use the power-convexity estimates in \cite[pp.~97--100]{Lindqvist2019} and set
\[
 \mathcal B_p(\xi,\zeta)
 =|\xi|^p-|\zeta|^p-p|\zeta|^{p-2}\zeta\cdot(\xi-\zeta).
\]
Testing with $\psi^pu$ gives
\[
 F_\psi=\int\psi^p|\nabla u|^p
 +p\int\psi^{p-1}u|\nabla u|^{p-2}\nabla u\cdot\nabla\psi.
\]
Since $\nabla v=\psi\nabla u+u\nabla\psi$,
\begin{equation}\label{eq-bregman-cutoff}
 R_\psi=E_v-F_\psi
 =\int\mathcal B_p(\nabla v,\psi\nabla u)\dd x\ge0.
\end{equation}

Testing with $v-z$ gives
\begin{equation}\label{eq-bregman-potential}
\begin{aligned}
 0\le\int\mathcal B_p(\nabla v,\nabla z)
 &=E_v-pF_\psi+(p-1)E_z\\
 &=R_\psi-(p-1)(F_\psi-E_z)\le R_\psi.
\end{aligned}
\end{equation}
For $p\ge2$, $\mathcal B_p(\xi,\zeta)\ge c_p|\xi-\zeta|^p$ gives \eqref{eq-h-energy}. For $1<p<2$,
\[
 \int\frac{|\nabla h|^2}{(|\nabla v|+|\nabla z|)^{2-p}}\dd x\le CR_\psi.
\]
The quotient is set to zero where $\nabla v=\nabla z=0$. H\"older yields
\[
 \int|\nabla h|^p
 \le CR_\psi^{p/2}
 \left(\int(|\nabla v|+|\nabla z|)^p\right)^{(2-p)/2}.
\]
Now $E_v+E_z\le2F_\psi+R_\psi$ and $\int V|h|^p\le\int|\nabla h|^p$.
\end{proof}

\begin{remark}
For $p=2$, the comparison estimate reduces to
\[
 \int_{\mathbb R^N}|\nabla(z-\psi u)|^2\dd x
 \le R_\psi
 =\int_{\mathbb R^N}u^2|\nabla\psi|^2\dd x.
\]
Thus the global gradient error between $z$ and $\psi u$
is controlled by the cutoff energy in the region
where $\nabla\psi\ne0$.
\end{remark}

\section{Directional pressure and source layers}\label{sec-pressure}

To transfer directional pressure from $z$ to $\psi u$, we combine Lemma~\ref{lem-potential} with source-layer endpoint estimates.

Write
\[
 \av_e H(e)=\frac1{|\mathbb S^{N-1}|}\int_{\mathbb S^{N-1}}H(e)\dd e.
\]

\subsection{Pressure from continuous Steiner symmetrization}

For $e\in\mathbb S^{N-1}$, write $x=y+te$, $y\in e^\perp$, and $\partial_e=e\cdot\nabla$. For $0\le g\in C_c^1$ and $w\in\Lip_{\mathrm{loc}}$, define
\begin{equation}\label{eq-pressure-def}
\begin{aligned}
S_e(g,w)=\int_{e^\perp}\dd y\int_0^\infty\dd\ell\sum_I
\biggl[&\int_I w(y+te)\dd t\\[-2pt]
&-\frac{|I|}{2}\bigl(w(y+a_Ie)+w(y+b_Ie)\bigr)\biggr].
\end{aligned}
\end{equation}
Here $I=(a_I,b_I)$ ranges over components of $\{t:g(y+te)>\ell\}$. Almost every level has finitely many regular endpoints. Using the Lipschitz constant $\Lip(w;B_{R_0})$, for $\supp g\subset B_{R_0}$,
\[
 \left|\int_Iw(y+te)\dd t-\frac{|I|}{2}\bigl(w(y+a_Ie)+w(y+b_Ie)\bigr)\right|
 \le CR_0\Lip(w;B_{R_0})|I|.
\]
Thus, by layer cake,
\[
\begin{aligned}
 |S_e(g,w)|
 &\le CR_0\Lip(w;B_{R_0})
 \int_{e^\perp}\dd y\int_0^\infty
 \sum_I |I|\dd\ell\\
 &=CR_0\Lip(w;B_{R_0})
 \int_{e^\perp}\dd y\int_{\mathbb R}\dd t
 \int_0^\infty
 \mathbf1_{\{g(y+te)>\ell\}}\dd\ell\\
 &=CR_0\Lip(w;B_{R_0})
 \int_{e^\perp}\dd y\int_{\mathbb R}g(y+te)\dd t\\
 &=CR_0\Lip(w;B_{R_0})
 \int_{\mathbb R^N}g(x)\dd x<\infty.
\end{aligned}
\]
Set
\begin{equation}\label{eq-stress}
 Q_e(\xi)=\frac{|\xi|^p}{p}-|\xi|^{p-2}(\xi\cdot e)^2,
 \qquad Q_e(0)=0.
\end{equation}

\begin{lemma}\label{lem-pressure}
If $0\le g\in C_c^1(\R^N)$ and $z=T_pg$, then
\begin{equation}\label{eq-pressure}
 S_e(g,z)\le\int_{\R^N}Q_e(\nabla z)\dd x.
\end{equation}
\end{lemma}

\begin{proof}
Since $g\in C_c^1$, local boundedness \cite[Theorem~1]{Serrin64} and gradient regularity \cite[Theorem~1]{Tolksdorf} give $z\in C^{1,\alpha}_{\mathrm{loc}}(\R^N)$ for some $\alpha\in(0,1)$.
We prove, by source and affine variations,
\[
 S_e(g,z)\le-\int g(x)(x\cdot e)\partial_ez
 =\int Q_e(\nabla z).
\]

Fix $e\in\mathbb S^{N-1}$ and $R_0>0$ with $\supp g\subset B_{R_0}$. For each $y\in e^\perp$ and level $\ell$, continuously Steiner-symmetrize $\{t:g(y+te)>\ell\}$ as in \cite[Section~2]{Brock}. Each interval $I=(a_I,b_I)$ has center $c_I=(a_I+b_I)/2$.
Here $\tau\ge0$ denotes the symmetrization parameter,
with $I(0)=I$ representing the original interval.
Until it meets another interval,
\[
 c_I(\tau)=e^{-\tau}c_I(0),\qquad |I(\tau)|=|I(0)|.
\]
The merging rule preserves length and nesting of the level sets. If $E_{\ell,\tau}(y)$ is the transformed set, define
\[
 g_\tau(y+te)=\int_0^\infty\mathbf1_{E_{\ell,\tau}(y)}(t)\dd\ell.
\]
For a nonnegative compactly supported Lipschitz function $\phi$, define $\phi_\tau$ in the same way. The rearrangement inequalities \cite{Brock1995,Brock2000,Brock} give
\[
 \int g_\tau\phi_\tau\ge\int g\phi,
 \qquad \int|\nabla\phi_\tau|^p\le\int|\nabla\phi|^p.
\]
Also $\supp g_\tau\subset B_{R_0}$ and $\norm{g_\tau-g}_\infty=O(\tau)$.

By \eqref{eq-intro-dual-functional}, for $z=T_pg$,
\begin{equation}\label{eq-dual-energy}
 W_p(g)=-J_g(z)=\frac{p-1}{p}\int|\nabla z|^p.
\end{equation}
For every such $\phi$,
\[
 \int g_\tau\phi_\tau-\frac1p\int|\nabla\phi_\tau|^p
 \ge\int g\phi-\frac1p\int|\nabla\phi|^p.
\]
Taking the supremum over $\phi$, using $|\phi|$ and density in $\D$, gives
\begin{equation}\label{eq-dual-monotonicity}
 W_p(g_\tau)\ge W_p(g).
\end{equation}

Hence, we prove the monotonicity of the dual energy.

For $z_\tau=T_pg_\tau$, maximality gives
\begin{equation}\label{eq-dual-difference-quotient}
 \int(g_\tau-g)z
 \le W_p(g_\tau)-W_p(g)
 \le\int(g_\tau-g)z_\tau.
\end{equation}
With $\D^*$ denoting the dual of $\D$, the common support gives
\[
 \norm{g_\tau-g}_{\D^*}=O(\tau),\qquad
 \norm{z_\tau}_{\D}^{p-1}\le\norm{g_\tau}_{\D^*}\le C.
\]
Testing the difference equation, with $A_p(\xi)=|\xi|^{p-2}\xi$,
\[
 \int\bigl(A_p(\nabla z_\tau)-A_p(\nabla z)\bigr)
       \cdot\nabla(z_\tau-z)
 =\int(g_\tau-g)(z_\tau-z)\longrightarrow0.
\]
The monotonicity inequalities \cite[Lemma~1.14.1]{zzt} give $z_\tau\to z$ in $\D$, hence
\[
 \left|\int(g_\tau-g)(z_\tau-z)\right|
 \le\norm{g_\tau-g}_{\D^*}\norm{z_\tau-z}_{\D}=o(\tau).
\]

At regular endpoints at $\tau=0$,
\[
 a_I'(0)=b_I'(0)=-c_I.
\]
Consequently,
\begin{equation}\label{eq-steiner-derivative}
 0\le\left.\frac{\dd}{\dd\tau}W_p(g_\tau)\right|_{0+}
 =-\int\dd y\dd\ell\sum_I c_I
 \bigl(z(y+b_Ie)-z(y+a_Ie)\bigr).
\end{equation}
To justify differentiation under the layer integral, write $z_y(t)=z(y+te)$. Interval transport preserves mass through mergers, with displacement at most $CR_0\tau$. Hence
\[
\begin{gathered}
 \frac1\tau\left|\int_{E_{\ell,\tau}(y)}z_y(t)\dd t
                    -\int_{E_{\ell,0}(y)}z_y(t)\dd t\right|
 \le CR_0\Lip(z;B_{R_0})|E_{\ell,0}(y)|,\\
 \int\dd y\dd\ell\,|E_{\ell,0}(y)|=\int g<\infty.
\end{gathered}
\]
Dominated convergence and \eqref{eq-dual-monotonicity} justify \eqref{eq-steiner-derivative}.

Integration by parts gives
\[
 \int_{a_I}^{b_I}z_y(t)\dd t
 -\frac{b_I-a_I}{2}\bigl(z_y(a_I)+z_y(b_I)\bigr)
 =-\int_{a_I}^{b_I}t z_y'(t)\dd t
 +c_I\bigl(z_y(b_I)-z_y(a_I)\bigr).
\]
Integrating over layers and using \eqref{eq-steiner-derivative},
\begin{equation}\label{eq-pressure-source-bound}
 S_e(g,z)\le-\int g(x)(x\cdot e)\partial_ez(x)\dd x.
\end{equation}
For $|t|<1/2$, set $z^{(t)}(x)=z(x+t(x\cdot e)e)$. Then
\[
 \left.\frac{\dd}{\dd t}\frac1p\int|\nabla z^{(t)}|^p\dd x\right|_{t=0}
 =\int\left(|\nabla z|^{p-2}(\partial_ez)^2-\frac{|\nabla z|^p}{p}\right)\dd x.
\]
Since $z^{(t)}\in\D$, $g\in C_c^1$, and $z\in C^1$, minimality gives
\[
 0=\left.\frac{\dd}{\dd t}J_g(z^{(t)})\right|_{t=0}
 =-\int Q_e(\nabla z)\dd x
  -\int g(x)(x\cdot e)\partial_ez\dd x.
\]
Together with \eqref{eq-pressure-source-bound}, this proves \eqref{eq-pressure}.
\end{proof}

\subsection{Length estimates with direction-dependent levels}

We use \eqref{eq-source-ball} to estimate long intervals, allowing the cutoff and endpoint level to depend on $e$.

Fix $s>0$ and a common annulus $\mathcal A\Subset B_s\setminus\{0\}$. For a measurable family of annular cutoffs, put
\[
 0\le\psi_e\le K,\qquad K>0,\qquad
 \supp\psi_e\subset\mathcal A,\qquad g_e=\psi_e^{p-1}f.
\]
For any $g\ge0$, let $L_g(e,x,t)\in[0,\infty]$ be the length of the component of $\{g>t\}\cap(x+\R e)$ containing $x$; set $L_g(e,x,t)=0$ if $g(x)\le t$.
Fix $x\ne0$, so $f(x),V(x)>0$.

For $\psi_e(x)>0$ and $\partial_eg_e(x)\ne0$, define the endpoint length by
\[
\begin{gathered}
 \varepsilon_e(x)=\operatorname{sgn}(\partial_eg_e(x)),\\
 L_e(x)=\sup\{L>0:g_e(x+t\varepsilon_e(x)e)>g_e(x)\text{ for }0<t<L\}.
\end{gathered}
\]
Otherwise set $L_e(x)=0$.

\begin{lemma}\label{lem-layer-tail}
For $\ell>0$ and fixed $0<\sigma<p$,
\begin{align}
 \av_e\bigl[\psi_e(x)^{p-1}\mathbf1_{\{L_e(x)>\ell\}}\bigr]
 &\le CK^{p-1}\min\{1,V(x)^{-1}\ell^{-p}\},\label{eq-layer-tail}\\
 \av_e\bigl[\psi_e(x)^{p-1}L_e(x)^\sigma\bigr]
 &\le C_\sigma K^{p-1}V(x)^{-\sigma/p}.\label{eq-layer-moment}
\end{align}
For $0<t<f(x)$,
\begin{equation}\label{eq-full-layer}
 \av_e L_f(e,x,t)\le Cu(x)^{(p-1)/p}t^{-1/p}.
\end{equation}
\end{lemma}
\begin{proof}
For $L_e(x)>\ell$,
\[
 f(x+s\varepsilon_e(x)e)>
 K^{1-p}\psi_e(x)^{p-1}f(x),\qquad 0<s<\ell.
\]
The equality $\varepsilon_e(x)e=\omega\in\mathbb S^{N-1}$
is possible only for $e=\omega$ or $e=-\omega$.
Hence each direction is counted at most twice. Polar integration gives
\[
 \int_{B_\ell(x)}f\dd y
 \ge cK^{1-p}f(x)\ell^N
 \av_e\bigl[\psi_e(x)^{p-1}\mathbf1_{\{L_e(x)>\ell\}}\bigr].
\]
This and \eqref{eq-source-ball} give \eqref{eq-layer-tail}. For $0<\sigma<p$,
\[
\begin{aligned}
 \av_e\bigl[\psi_e^{p-1}L_e^\sigma\bigr]
 &=\sigma\int_0^\infty\ell^{\sigma-1}
   \av_e\bigl[\psi_e^{p-1}\mathbf1_{\{L_e>\ell\}}\bigr]\dd\ell\\
 &\le CK^{p-1}\left(
 \int_0^{V^{-1/p}}\ell^{\sigma-1}\dd\ell
 +V^{-1}\int_{V^{-1/p}}^\infty\ell^{\sigma-p-1}\dd\ell\right)\\
 &\le C_\sigma K^{p-1}V^{-\sigma/p}.
\end{aligned}
\]

Since $L_f>\ell$ implies one side has length $>\ell/2$, the same argument gives
\[
 \av_e\mathbf1_{\{L_f(e,x,t)>\ell\}}
 \le C\min\{1,u(x)^{p-1}/(t\ell^p)\}.
\]
Integration in $\ell$ proves \eqref{eq-full-layer} since $p>1$.
\end{proof}

\subsection{Endpoint errors}

For each direction, write
\[
 z_e=T_pg_e,\qquad v_e=\psi_eu,\qquad h_e=z_e-v_e,
\]
with $F_{\psi_e},R_{\psi_e},E_{z_e},E_{v_e}$ defined by \eqref{eq-potential-notation} and Lemma~\ref{lem-potential} with $\psi=\psi_e$. Fix
\begin{equation}\label{eq-sigma}
 p'=\frac p{p-1},\qquad 1<\sigma<\min\{p,p'\},\qquad
 \sigma'=\frac\sigma{\sigma-1}>p.
\end{equation}
The length moment uses $\sigma<p$; interpolation uses $\sigma'>p$. These exponents are fixed for all cutoffs.

We write $C_{\psi}$ for constants that depend polynomially on $\sup_e\norm{\psi_e^t}_{W^{1,\infty}}$ for finitely many fixed $t>0$, and may increase from line to line. They are independent of $u$; their dependence on the cutoff parameters is estimated in \eqref{eq-P}.

Layer cake and one-dimensional coarea give

\begin{align}
 S_e(g_e,z_e)-S_e(g_e,v_e)
 &=\int_{e^\perp}\dd y\int_0^\infty\dd\ell
 \sum_I\biggl[\int_I h_e(y+te)\dd t-\frac{|I|}{2}
 \bigl(h_e(y+a_Ie)+h_e(y+b_Ie)\bigr)\biggr]\notag\\
 &=\int g_eh_e\dd x
 -\frac12\int h_e|\partial_eg_e|L_e\dd x.
\end{align}

Hence
\begin{align}
 |S_e(g_e,z_e)-S_e(g_e,v_e)|
 &\le\int g_e|h_e|\dd x+\nu_e,\label{eq-endpoint-difference}\\
 \nu_e&=\frac12\int|h_e|\,|\partial_eg_e|L_e\dd x.\label{eq-endpoint-def}
\end{align}

On $\mathcal A$,
\begin{equation}\label{eq-source-gradient}
\begin{aligned}
 \nabla f&=f\left(a\frac{x}{|x|^2}+q\frac{\nabla u}{u}\right),
 &|\nabla f|&\le C\frac fu(|\nabla u|+u),\\
 \nabla g_e&=\psi_e^{p-1}\nabla f+f\nabla(\psi_e^{p-1}).
\end{aligned}
\end{equation}
Here $C=C(\mathcal A,a,q)$, and $\nabla(\psi_e^{p-1})=0$ on $\{\psi_e=0\}$.

Factor
\[
 |\partial_eg_e|
 \le\psi_e^{(p-1)/\sigma}
 \left[\psi_e^{(p-1)/\sigma'}|\partial_ef|
       +\sigma'f\left|\partial_e\bigl(\psi_e^{(p-1)/\sigma'}\bigr)\right|\right].
\]
Directional H\"older, \eqref{eq-layer-moment}, and \eqref{eq-source-gradient} yield
\begin{equation}\label{eq-endpoint-pointwise-first}
 \av_e\bigl[|h_e|\,|\partial_eg_e|L_e\bigr]
 \le C_{\psi}\frac fu(|\nabla u|+u)V^{-1/p}
       (\av_e|h_e|^{\sigma'})^{1/\sigma'}.
\end{equation}

Since $|h_e|\le Ku$ and $\sigma'>p$,
\[
 (\av_e|h_e|^{\sigma'})^{1/\sigma'}
 \le(Ku)^{1-p/\sigma'}(\av_e|h_e|^p)^{1/\sigma'}.
\]
Since $V=f/u^{p-1}$,
\[
 \frac fu\,V^{-1/p}u^{1-p/\sigma'}V^{-1/\sigma'}
 =f^{1-1/p-1/\sigma'}u^{(p-1)/p-1/\sigma'}
 =(fu)^{1/\sigma-1/p}.
\]
Thus \eqref{eq-endpoint-pointwise-first} is bounded by
\[
 C_{\psi}\bigl(V\av_e|h_e|^p\bigr)^{1/\sigma'}
   (fu)^{1/\sigma-1/p}(|\nabla u|+u).
\]
Now
\[
 \frac1\sigma-\frac1p>0,\qquad
 \frac1{\sigma'}+\left(\frac1\sigma-\frac1p\right)+\frac1p=1.
\]
Spatial H\"older and Fubini give
\begin{equation}\label{eq-endpoint-bound}
 \av_e\nu_e\le C_{\psi}
 \left(\av_e\int V|h_e|^p\right)^{1/\sigma'}
 F(s)^{1/\sigma-1/p}
 \left(\int_{\mathcal A}(|\nabla u|^p+u^p)\right)^{1/p}.
\end{equation}

\subsection{Recovery of the reaction term}

We recover the coefficient $1/(q+1)$ in Pohozaev's identity.

Since $d=q-p+1>0$,
\begin{equation}\label{eq-source-inversion}
 v_e=b_eg_e^{1/q},\qquad b_e=\psi_e^{d/q}|x|^{-a/q}.
\end{equation}
The zero extension satisfies $b_e\in\Lip(\R^N)$. On layer endpoints, $v_e=b_e\ell^{1/q}$.

Layer cake gives
\begin{align*}
 \int\dd y\dd\ell\sum_I\ell^{1/q}\int_I b_e(y+te)\dd t
 &=\int b_e(x)\int_0^{g_e(x)}\ell^{1/q}\dd\ell\dd x\\
 &=\frac q{q+1}\int b_eg_e^{1+1/q}\dd x
 =\frac q{q+1}\int g_ev_e\dd x.
\end{align*}
Therefore
\begin{equation}\label{eq-source-recovery}
 S_e(g_e,v_e)=\frac1m\int g_ev_e\dd x+\mathcal E_e,
\end{equation}
where, with $x=y+te$,
\begin{equation}\label{eq-recovery-error-definition}
 \mathcal E_e=\int\dd y\dd\ell\sum_I\ell^{1/q}
 \int_I\left[b_e(x)-\frac{b_e(y+a_Ie)+b_e(y+b_Ie)}2\right]\dd t.
\end{equation}
By Lipschitz continuity,
\[
 |\mathcal E_e|\le\Lip(b_e)
 \int_{\mathcal A}\dd x\int_0^{g_e(x)}\ell^{1/q}L_{g_e}(e,x,\ell)\dd\ell.
\]
Since $g_e\le K^{p-1}f$,
\[
 L_{g_e}(e,x,\ell)\le L_f(e,x,\ell/K^{p-1}).
\]
Enlarging the upper limit and using \eqref{eq-full-layer},
\begin{align*}
 \av_e|\mathcal E_e|
 &\le C_{\psi}\int_{\mathcal A}u^{(p-1)/p}
       \int_0^{K^{p-1}f}\ell^{1/q-1/p}\dd\ell\dd x\\
 &\le C_{\psi}\int_{\mathcal A}u^{(p-1)/p}f^{1+1/q-1/p}\dd x.
\end{align*}
On $\mathcal A$,
\[
\begin{gathered}
 1+\frac1q-\frac1p>0,\qquad 0<c\le |x|^a\le C,\\
 \frac{p-1}{p}+q\left(1+\frac1q-\frac1p\right)
 =m-\frac dp=m(1-d_0)<m.
\end{gathered}
\]
H\"older gives
\begin{equation}\label{eq-recovery-error}
 \av_e|\mathcal E_e|\le C_{\psi} F(s)^{1-d_0}.
\end{equation}

\section{Local energy closure and the scalar theorem}\label{sec-close}

We combine averaged directional pressure with Pohozaev's identity, then absorb the localization errors.

\subsection{Annular cutoffs and directional concentration}

Fix $1\le r<s\le2$, $\delta=s-r$, and a smooth nonincreasing radial cutoff satisfying
\[
 \eta=1\text{ on }B_r,\qquad \eta=0\text{ outside }B_s,\qquad
 \supp\eta'\subset\{r+\delta/4<|x|<r+\delta/2\}.
\]
Set $w(\rho)=-\rho\eta'(\rho)\ge0$. For $e,\omega\in\mathbb S^{N-1}$ and $0<\kappa<1/4$, take
\[
 k_\kappa(e,\omega)=
 \frac{\vartheta(|e-\omega|^2/\kappa^2)}
 {\av_{\tilde e}\vartheta(|\tilde e-\omega|^2/\kappa^2)},
\]
where
\[
 \vartheta(t)=\begin{cases}e^{-1/(1-t)},&0\le t<1,\\0,&t\ge1.\end{cases}
\]
The denominator is $\asymp\kappa^{N-1}$, and
\[
 \supp_e k_\kappa\subset\{|e-\omega|\le\kappa\},\qquad
 \av_e k_\kappa(e,\omega)=1.
\]
Using bumps with smooth positive powers, set
\begin{equation}\label{eq-cutoff-family}
 \psi_e(x)^p=w(|x|)k_\kappa(e,x/|x|).
\end{equation}
For $\omega=x/|x|$ and $\xi\in\R^N$, the support and normalization of $k_\kappa$ give
\[
 \av_e\psi_e^p=w,\qquad
 \left|\av_e\psi_e^pQ_e(\xi)-wQ_\omega(\xi)\right|
 \le C\kappa w|\xi|^p.
\]

For fixed $t>0$,
\begin{align}
 \norm{\psi_e^t}_\infty
 &\le C_t\delta^{-t/p}\kappa^{-(N-1)t/p},\label{eq-cutoff-size}\\
 \norm{\nabla(\psi_e^t)}_\infty
 &\le C_t\delta^{-t/p}\kappa^{-(N-1)t/p}
 (\delta^{-1}+\kappa^{-1}).\label{eq-cutoff-derivative}
\end{align}
The cutoff bounds \eqref{eq-cutoff-size}--\eqref{eq-cutoff-derivative} imply
\begin{equation}\label{eq-P}
 C_{\psi}\le C(\delta\kappa)^{-M},\qquad 0<M<\infty,
\end{equation}
with $C,M$ independent of $u,\delta,\kappa$ and subsequent absorption parameters.

Set $X=1+F(s)$. Choose $\mathcal A$ containing all cutoff supports, with a $c\delta$ buffer inside $B_s\setminus\overline{B_1}$. Then
\begin{equation}\label{eq-annular-energy}
 \int_{\mathcal A}u^p\le CF(s)^{p/m},\qquad
 \int_{\mathcal A}|\nabla u|^p\le C\delta^{-p}X.
\end{equation}
Since $0<c\le|x|^a\le C$ on $\mathcal A$,
\[
 \int_{\mathcal A}u^p
 \le |\mathcal A|^{1-p/m}\left(\int_{\mathcal A}u^m\right)^{p/m}
 \le CF(s)^{p/m}.
\]
Choose $\xi\in C_c^\infty(B_s\setminus\overline{B_1})$, $\xi=1$ on $\mathcal A$, $|\nabla\xi|\le C/\delta$. Testing with $\xi^pu$ gives
\[
 \int\xi^p|\nabla u|^p
 \le F(s)+\frac12\int\xi^p|\nabla u|^p
       +C\delta^{-p}\int_{\supp\xi}u^p.
\]
Since $\delta\le1$ and $p/m<1$,
\[
 \int_{\mathcal A}|\nabla u|^p
 \le CF(s)+C\delta^{-p}F(s)^{p/m}
 \le C\delta^{-p}(1+F(s)).
\]

\subsection{Fixed sublinear error exponents}

All error bounds below use one fixed exponent $\alpha_*<1$. Enlarge $C,M$ in \eqref{eq-P} to cover all cutoff constants.

For $p\ge2$, \eqref{eq-bregman-cutoff} gives
\[
 R_{\psi_e}\le C\int\left(
 |\psi_e\nabla u|^{p-2}|u\nabla\psi_e|^2+|u\nabla\psi_e|^p\right).
\]
By \eqref{eq-annular-energy},
\[
 R_{\psi_e}\le C_{\psi}\bigl(X^{(p-2)/p+2/m}+X^{p/m}\bigr)
 =C_{\psi}\bigl(X^{1-2d_0}+X^{1-pd_0}\bigr)
 \le 2C_{\psi}X^{1-2d_0}\qquad(p\ge2).
\]
For $1<p<2$, $R_{\psi_e}\le C\int|u\nabla\psi_e|^p$. Hence
\begin{equation}\label{eq-cutoff-remainder-bound}
 R_{\psi_e}\le
 \begin{cases}
 C_{\psi}X^{1-2d_0},&p\ge2,\\
 C_{\psi}X^{1-pd_0},&1<p<2.
 \end{cases}
\end{equation}
Since $F_{\psi_e}+E_{z_e}+E_{v_e}\le C_{\psi}X$, Lemma~\ref{lem-potential} gives
\begin{equation}\label{eq-h-final}
 \av_e\int\bigl(|\nabla h_e|^p+V|h_e|^p\bigr)\dd x
 \le C_{\psi}X^{1-\varepsilon_0},\qquad
 \varepsilon_0=\begin{cases}2d_0,&p\ge2,\\p^2d_0/2,&1<p<2.\end{cases}
\end{equation}
For $1<p<2$,
\[
 \frac p2(1-pd_0)+\frac{2-p}{2}=1-\frac{p^2d_0}{2}.
\]

We also use the stress estimate
\begin{equation}\label{eq-stress-continuity}
 |Q_e(\xi)-Q_e(\zeta)|
 \le C_p(|\xi|+|\zeta|)^{p-1}|\xi-\zeta|.
\end{equation}
H\"older, $g_e=Vv_e^{p-1}$, \eqref{eq-h-final}, and $\nabla v_e-\psi_e\nabla u=u\nabla\psi_e$ give the following three groups of errors:
\begin{equation}\label{eq-three-error-groups}
\begin{aligned}
 \av_e\nu_e&\le C_{\psi}X^{1-\varepsilon_0/\sigma'},\\
 \av_e\left[\int g_e|h_e|
 +\int|Q_e(\nabla z_e)-Q_e(\nabla v_e)|\right]
 &\le C_{\psi}X^{1-\varepsilon_0/p},\\
 \av_e\left[\int|Q_e(\nabla v_e)-\psi_e^pQ_e(\nabla u)|
 +|\mathcal E_e|\right]
 &\le C_{\psi}X^{1-d_0}.
\end{aligned}
\end{equation}
The first line is \eqref{eq-endpoint-bound}; the remaining exponents are
\[
 \frac{p-1}{p}+\frac{1-\varepsilon_0}{p}=1-\frac{\varepsilon_0} p,
 \qquad \frac{p-1}{p}+\frac1m=1-d_0,
\]
with \eqref{eq-recovery-error} supplying the last term.
Fix
\begin{equation}\label{eq-alpha-star}
 \alpha_*=1-\frac{\varepsilon_0}{\sigma'},\qquad 0<\alpha_*<1.
\end{equation}
Indeed,
\[
 0<\varepsilon_0<1,\qquad \frac{\varepsilon_0}{\sigma'}<\frac{\varepsilon_0} p\le d_0,
 \qquad
 \max\{1-\varepsilon_0/\sigma',1-\varepsilon_0/p,1-d_0\}=\alpha_*.
\]
The exponents $\alpha_*,M$ are independent of $u,\delta,\kappa$.

Combine \eqref{eq-source-recovery}, \eqref{eq-endpoint-difference}, and \eqref{eq-pressure}:
\[
\begin{aligned}
 \frac1m\int g_ev_e
 &=S_e(g_e,v_e)-\mathcal E_e\\
 &\le\int Q_e(\nabla z_e)
       +\int g_e|h_e|+\nu_e+|\mathcal E_e|.
\end{aligned}
\]
By \eqref{eq-three-error-groups}, $g_ev_e=\psi_e^pfu$, and $\av_e\psi_e^p=w$,
\begin{equation}\label{eq-averaged-pressure}
 \frac1m\int wfu
 \le\av_e\int\psi_e^pQ_e(\nabla u)
 +C(\delta\kappa)^{-M}X^{\alpha_*}.
\end{equation}
Write $u_r=\nabla u\cdot x/|x|$. The stress concentration following \eqref{eq-cutoff-family}, $\norm w_\infty\le C/\delta$, and \eqref{eq-annular-energy} give
\begin{equation}\label{eq-radial-pressure}
\begin{aligned}
 \frac1m\int wfu
 &\le\int w\left(\frac{|\nabla u|^p}{p}-|\nabla u|^{p-2}u_r^2\right)\\
 &\quad+C\kappa\delta^{-p-1}X+C(\delta\kappa)^{-M}X^{\alpha_*},
\end{aligned}
\end{equation}
Gradient products are zero at $\nabla u=0$.

\subsection{Localized Pohozaev identity and removal of the origin}

The boundary combination in Pohozaev's identity is controlled by \eqref{eq-radial-pressure}.

Unlike $F_\psi$ in \eqref{eq-potential-notation}, the following integrals use the cutoffs linearly:
\[
\begin{gathered}
 F_\eta=\int\eta fu,\qquad E_\eta=\int\eta|\nabla u|^p,\\
 F_w=\int wfu,\qquad E_w=\int w|\nabla u|^p,\qquad
 T_w=\int w|\nabla u|^{p-2}u_r^2.
\end{gathered}
\]
With $DY=(\partial_jY_i)_{i,j=1}^N$ and $I=\mathrm{Id}_{\R^N}$,
\[
 D(x\eta)=\eta I-w\,\omega\otimes\omega,\qquad
 \diver(x\eta)=N\eta-w,\qquad
 x\eta\cdot\nabla(|x|^a)=a\eta|x|^a.
\]
With $s_+=\max\{s,0\}$, apply \cite[Lemma~1]{DMS} to
\[
 L(x,s,\xi)=\frac{|\xi|^p}{p}-\frac{|x|^a(s_+)^m}{m},\qquad x\ne0.
\]
It is $C^1$ and strictly convex in $\xi$; $u$ is locally $C^1$. For $Y\in C_c^1(\R^N\setminus\{0\};\R^N)$,
\[
 \int\left(\frac{\diver Y}{p}|\nabla u|^p
 -|\nabla u|^{p-2}\nabla u\cdot DY\nabla u\right)
 =\frac1m\int\left((\diver Y)|x|^a+Y\cdot\nabla(|x|^a)\right)u^m.
\]
Take $Y=x\eta(|x|)\chi(|x|/\eps)$, with $\chi=0$ on $[0,1]$ and $\chi=1$ on $[2,\infty)$. Since $|x\nabla\chi(|x|/\eps)|\le C$, the cutoff terms are bounded by
\[
 C\int_{B_{2\eps}}\bigl(|\nabla u|^p+|x|^au^m\bigr)\dd x\longrightarrow0.
\]
Letting $\eps\downarrow0$ gives the Localized Pohozaev identity:
\begin{equation}\label{eq-pohozaev}
 \frac{N+a}{m}F_\eta-\frac{N-p}{p}E_\eta
 =\frac1mF_w-\frac1pE_w+T_w.
\end{equation}

Testing with $\eta u$ gives
\[
 E_\eta=F_\eta-J_\eta,\qquad
 J_\eta:=\int u|\nabla u|^{p-2}\nabla u\cdot\nabla\eta.
\]
H\"older and \eqref{eq-annular-energy} imply
\[
 |J_\eta|\le C\delta^{-1}
 \left(\int_{\mathcal A}|\nabla u|^p\right)^{(p-1)/p}
 \left(\int_{\mathcal A}u^p\right)^{1/p}
 \le C_{\psi}X^{(p-1)/p+1/m}=C_{\psi}X^{1-d_0}.
\]
By \eqref{eq-pohozaev} and \eqref{eq-radial-pressure},
\[
\begin{aligned}
 \Gamma F(r)\le\Gamma F_\eta
 &=\frac{F_w}{m}-\frac{E_w}{p}+T_w-\frac{N-p}{p}J_\eta\\
 &\le C\kappa\delta^{-p-1}X
       +C(\delta\kappa)^{-M}X^{\alpha_*}+C|J_\eta|.
\end{aligned}
\]
Since $F_\eta\ge F(r)$ and $1-d_0\le\alpha_*$,
\begin{equation}\label{eq-core-energy}
\begin{aligned}
 \Gamma F(r)&\le C\kappa\delta^{-p-1}(1+F(s))\\
 &\quad+C(\delta\kappa)^{-M}(1+F(s))^{\alpha_*},
 \qquad 1\le r<s\le2.
\end{aligned}
\end{equation}
Here $C$ is independent of $u$.

\subsection{Hole filling and the scaling contradiction}

We absorb the linear and sublinear errors in \eqref{eq-core-energy}.

Absorb $\Gamma$ into $C$. For small $\theta>0$, choose
\[
 \kappa=c_0\theta\delta^{p+1},\qquad
 C\kappa\delta^{-p-1}\le\frac\theta4,\qquad
 C(\delta\kappa)^{-M}\le C\theta^{-M}\delta^{-M(p+2)},
\]
with fixed small $c_0>0$. Put $A_0=C\theta^{-M}\delta^{-M(p+2)}$. Young gives
\[
 A_0X^{\alpha_*}\le\frac\theta4 X
 +C\theta^{-\alpha_*/(1-\alpha_*)}A_0^{1/(1-\alpha_*)}.
\]
Hence
\begin{equation}\label{eq-hole}
 F(r)\le\theta F(s)+C\theta^{-A}(s-r)^{-B},\qquad
 A=\frac{M+\alpha_*}{1-\alpha_*},\quad
 B=\frac{M(p+2)}{1-\alpha_*}.
\end{equation}
Fix $\theta2^B<1$ and $r_j=2-2^{-j}$. For each integer $K\ge1$, iteration gives
\[
 F(1)\le\theta^KF(r_K)
 +C\theta^{-A}2^B\sum_{j=0}^{K-1}(\theta2^B)^j.
\]
Since $F(2)<\infty$,
\[
 0\le\theta^KF(r_K)\le\theta^KF(2)\longrightarrow0,
 \qquad
 \sum_{j=0}^\infty(\theta2^B)^j=\frac1{1-\theta2^B}<\infty.
\]
Letting $K\to\infty$ proves
\begin{equation}\label{eq-unit-energy}
 F(1)\le C_*.
\end{equation}
Here $C_*=C_*(N,p,a,q)$ is independent of $u$.

\begin{proof}[Proof of Theorem~\ref{thm-scalar}]
Put
\[
\begin{gathered}
 \beta=\frac{p+a}{q-p+1},\qquad (p-1)\beta+p+a=q\beta,\\
 \Lambda=\beta(q+1)-(N+a)=p\beta-(N-p)
 =\frac{(N-p)(q_H(a)-q)}{q-p+1}>0.
\end{gathered}
\]
The rescaling $u_R(x)=R^\beta u(Rx)$ preserves \eqref{eq-solution-class}--\eqref{eq-weak}. By \eqref{eq-unit-energy},
\[
 C_*\ge\int_{B_1}|x|^au_R(x)^{q+1}\dd x
 =R^\Lambda F(R).
\]
For $R\ge1$,
\[
 0<F(1)\le F(R)\le C_*R^{-\Lambda}\longrightarrow0,
\]
a contradiction.
\end{proof}

\subsection{The critical endpoint and the solution class}

At $q=q_H(a)$, $\Lambda=0$ and the scaling contradiction ends. With $C_0>0$ normalized below, the critical profile \cite[Eq.~(1.4)]{PS},
\begin{equation}\label{eq-bubble}
\begin{gathered}
 U(x)=C_0(1+|x|^{(p+a)/(p-1)})^{-(N-p)/(p+a)},\\
 C_0^{q_H(a)-p+1}=(N+a)\left(\frac{N-p}{p-1}\right)^{p-1}.
\end{gathered}
\end{equation}
satisfies the equation in both solution classes. For $r=|x|$, put
\[
 t=\frac{p+a}{p-1}>0,\qquad b=\frac{N-p}{p+a}>0,\qquad bt=\frac{N-p}{p-1}.
\]
Then
\[
 U'(r)=-C_0bt\,r^{t-1}(1+r^t)^{-b-1},
\]
and
\[
 -\Delta_pU=(C_0bt)^{p-1}(N+a)r^a
 (1+r^t)^{-(b+1)(p-1)-1}.
\]
Since $bq_H(a)=(b+1)(p-1)+1$, \eqref{eq-bubble} gives $-\Delta_pU=r^aU^{q_H(a)}$ for $r>0$.

Moreover, $U\in C^\infty(\R^N\setminus\{0\})\cap C(\R^N)$, $U(0)=C_0$, and
\[
\begin{aligned}
 \int_{B_\eps}|\nabla U|^p
 &=O\bigl(\eps^{N+p(t-1)}\bigr),& N+p(t-1)&=N-p+pt>0,\\
 \int_{\partial B_\eps}|\nabla U|^{p-1}
 &=O(\eps^{N+a}),& N+a&>0,\\
 \int_{B_\eps}|x|^aU^{q_H(a)}&=O(\eps^{N+a}).
\end{aligned}
\]
Hence $U\in W^{1,p}_{\mathrm{loc}}\cap L^\infty_{\mathrm{loc}}$ and satisfies \eqref{eq-weak}. This proves sharpness in both classes.

The origin regularity is
\[
\begin{aligned}
 U\in C^1(\R^N)&\quad\Longleftrightarrow\quad t>1
                  \quad\Longleftrightarrow\quad a>-1,\\
 U\in C^2(\R^N)&\quad\Longleftrightarrow\quad t\ge2
                  \quad\Longleftrightarrow\quad a\ge p-2,\\
 |\nabla U|^{p-2}\nabla U\in C^1(\R^N)
              &\quad\Longleftrightarrow\quad a\ge0.
\end{aligned}
\]

For $-N<a\le-p$, a positive solution in \eqref{eq-classical-class} would satisfy \eqref{eq-weak} by Section~\ref{sec-origin}. The Wolff bound gives
\[
 u(0)\ge c\int_0^\eps r^{(p+a)/(p-1)}\frac{\dd r}{r}=\infty.
\]
For $a\le-N$, $u(0)>0$ implies $|x|^au^q\notin L^1_{\mathrm{loc}}$.

\section{Critical classification for negative weights}\label{sec-negative-critical}

The notation in this section is local to the critical argument. Here $g$ denotes the cone metric in \eqref{eq-cone-data}; in Sections~\ref{sec-basic}--\ref{sec-close} it denoted the auxiliary source. Assume
\begin{equation}\label{eq-critical-parameters}
\begin{gathered}
 N\ge3,\quad 1<p<N,\quad -p<a<0,\qquad
 m=\frac{p(N+a)}{N-p},\quad p^*=\frac{Np}{N-p},\\
 \beta=\frac{N-p}{p},\qquad h_a=\frac{p+a}{p}\in(0,1),\qquad
 -\Delta_pu=|x|^au^{m-1}.
\end{gathered}
\end{equation}
Throughout, $u>0$, $u\in W^{1,p}_{\mathrm{loc}}(\R^N)$, and the weak identity holds on all of $\R^N$.

\subsection{Estimates derived from the weak equation}

\begin{lemma}\label{lem-critical-prior}
For some $\alpha>0$, the solution has a representative satisfying
\[
 u\in C(\R^N)\cap L^\infty_{\mathrm{loc}}(\R^N)
       \cap C^{1,\alpha}_{\mathrm{loc}}(\R^N\setminus\{0\}),
 \qquad u>0\text{ in }\R^N.
\]
For $r=|x|\ge1$,
\begin{equation}\label{eq-critical-coarse}
 c_u r^{-(N-p)/(p-1)}\le u(x)\le C_u r^{-\beta},
 \qquad |\nabla u(x)|\le C_u\frac{u(x)}r.
\end{equation}
If
\begin{equation}\label{eq-apex-exponent}
 h_a<\alpha_1<\min\left\{1,\frac{p+a}{p-1}\right\},
\end{equation}
then, near the origin,
\begin{equation}\label{eq-critical-apex}
 |u(x)-u(0)|\le C_u r^{\alpha_1},\qquad
 |\nabla u(x)|\le C_u\left(r^{\alpha_1-1}+r^{(a+1)/(p-1)}\right).
\end{equation}
\end{lemma}

\begin{proof}

Put $\sigma=-a\in(0,p)$ and $\omega=|x|^{-\sigma}$. Hardy and Sobolev, followed by H\"older, give
\begin{equation}\label{eq-critical-HS}
 \left(\int\omega|z|^m\right)^{p/m}
 \le C\int|\nabla z|^p,
 \qquad z\in W^{1,p}(\R^N),\quad \supp z\text{ compact}.
\end{equation}
Taking $z=\chi u$ with $\chi\in C_c^\infty$ yields $\omega u^m\in L^1_{\mathrm{loc}}$. For $\phi\in C_c^\infty$,
\[
 \left|\int\omega u^{m-1}\phi\right|
 \le\left(\int_{\supp\phi}\omega u^m\right)^{(m-1)/m}
      \left(\int\omega|\phi|^m\right)^{1/m}
 \le C\norm{\nabla\phi}_{L^p};
\]
the weak identity extends to compactly supported $W^{1,p}$ tests.

Fix $t\ge1$, $L,R>0$, $u_L=\min\{u,L\}$, and $0\le\chi\in C_c^\infty(B_R)$. Test with $\chi^puu_L^{p(t-1)}$ and set $z_L=\chi uu_L^{t-1}$. Product differentiation and Young's inequality give
\begin{equation}\label{eq-critical-bootstrap}
 \int|\nabla z_L|^p
 \le C_{p,t}\left(\int u^{pt}|\nabla\chi|^p
                 +\int\omega u^{m-p}z_L^p\right).
\end{equation}
The test is admissible since
\[
 \int\omega u^{m-p}z_L^p
 =\int\chi^p\omega u^m u_L^{p(t-1)}
 \le L^{p(t-1)}\int_{B_R}\omega u^m<\infty.
\]
On $\{u>K\}$,
\begin{align*}
 \int_{\{u>K\}}\omega u^{m-p}z_L^p
 &\le\left(\int_{B_R\cap\{u>K\}}\omega u^m\right)^{(m-p)/m}
       \left(\int\omega z_L^m\right)^{p/m}\\
 &\le C\left(\int_{B_R\cap\{u>K\}}\omega u^m\right)^{(m-p)/m}
       \int|\nabla z_L|^p.
\end{align*}
For fixed $t$, choose $K$ large enough that
\[
 C_{p,t}\left(\int_{B_R\cap\{u>K\}}\omega u^m\right)^{(m-p)/m}
 \le\frac12.
\]
On $\{u\le K\}$, use $u^{m-p}\le K^{m-p}$ and
\begin{equation}\label{eq-critical-small-weight}
\begin{aligned}
 \int|x|^{-\sigma}|z|^p
 &\le\left(\int|x|^{-p}|z|^p\right)^{\sigma/p}
       \left(\int|z|^p\right)^{1-\sigma/p}\\
 &\le\eps\int|\nabla z|^p+C_\eps\int|z|^p.
\end{aligned}
\end{equation}
Consequently,
\[
 \int|\nabla z_L|^p
 \le C_{p,t,u,R}\int_{B_R}u^{pt}(|\nabla\chi|^p+\chi^p),
\]
uniformly in $L$. Sobolev and Fatou imply
\[
 u\in L^{pt}(B_R)\quad\Longrightarrow\quad
 u\in L^{p^*t}(B_r),\qquad B_r\Subset B_R.
\]
Iteration on nested balls gives
\[
 t_j=\left(\frac{p^*}{p}\right)^j\longrightarrow\infty,\qquad
 u\in\bigcap_{j\ge0}L^{pt_j}_{\mathrm{loc}}
 \quad\Longrightarrow\quad
 u\in\bigcap_{1\le\ell<\infty}L^\ell_{\mathrm{loc}}.
\]

Choose $N/p<s<t_0<N/\sigma$. Since $\omega\in L^{t_0}_{\mathrm{loc}}$, H\"older gives $f=\omega u^{m-1}\in L^s_{\mathrm{loc}}$, hence
\[
 \int_{B_\rho(x)}f\le C\rho^{N-N/s}
 \le C\rho^{N-p+\alpha_0(p-1)},\qquad
 0<\alpha_0<\min\left\{1,\frac{p-N/s}{p-1}\right\}.
\]
The local measure-growth criterion of Kilpel\"ainen--Zhong \cite{KZ2002}, in the form \cite[Theorem~2.9]{Kilpelainen}, yields $u\in C^{0,\alpha_0}_{\mathrm{loc}}$. Here the homogeneous H\"older exponent in that criterion is $1$ for the $p$-Laplacian. The continuous representative is positive by weak Harnack. Interior regularity \cite{DiBenedetto1983,Tolksdorf} gives $C^{1,\alpha}$ regularity away from zero.

Here $a<0$ implies $p<m<p^*$. Set
\[
 u_R(y)=R^\beta u(Ry),\qquad
 -\Delta_pu_R=|y|^au_R^{m-1}.
\]
If their suprema on $\{1/2\le|y|\le2\}$ were unbounded, apply \cite[Lemma~5.1]{PQS2007} to $M=u_R^{(m-p)/p}$ on $\{1/4<|y|<4\}$. At the selected points $y_j$, put $K_j=M(y_j)$ and
\[
 z_j(\xi)=\frac{u_{R_j}(y_j+K_j^{-1}\xi)}{u_{R_j}(y_j)}.
\]
The doubling lemma and local compactness give, along a subsequence,
\[
\begin{gathered}
 z_j(0)=1,\qquad 0\le z_j\le C\text{ in }B_j,\qquad
 z_j\longrightarrow z\text{ in }C^1_{\mathrm{loc}},\\
 z(0)=1,\qquad 0\le z\le C,\qquad
 -\Delta_pz=c z^{m-1}\text{ in }\R^N,\quad c>0.
\end{gathered}
\]
By \cite[Theorem~III]{SZ}, $z\equiv0$, contradicting $z(0)=1$. Hence
\[
 \sup_{R\ge1}\sup_{1/2\le|y|\le2}u_R(y)\le C_u
 \quad\Longrightarrow\quad
 u(x)\le C_u|x|^{-\beta}\quad(|x|\ge1).
\]

The potential in $-\Delta_pu=Vu^{p-1}$ satisfies
\[
 V=|x|^au^{m-p}\le C_ur^{-p},\qquad \beta(m-p)=p+a.
\]
Annular Harnack and the interior gradient estimate give $|\nabla u(x)|\le C_uu(x)/r$. Comparison with the exterior $p$-harmonic function $c_ur^{-(N-p)/(p-1)}$ gives the lower bound in \eqref{eq-critical-coarse}.

By local boundedness and \eqref{eq-apex-exponent}, uniformly for small balls,
\[
 \alpha_1(p-1)<p+a,\qquad
 \int_{B_\rho(x)}f\le C\rho^{N+a}
 \le C\rho^{N-p+\alpha_1(p-1)}.
\]
The measure-growth criterion gives $u\in C^{0,\alpha_1}$ near zero. The interior gradient estimate on $B_{r/2}(x)$, $r=|x|$, gives
\[
\begin{aligned}
 |\nabla u(x)|&\le C\left(
 r^{-1}\operatorname{osc}_{B_{r/2}(x)}u
 +r^{1/(p-1)}\norm f_{L^\infty(B_{r/2}(x))}^{1/(p-1)}\right)\\
 &\le C\left(r^{\alpha_1-1}+r^{(a+1)/(p-1)}\right).
\end{aligned}
\]
\end{proof}

\subsection{The weighted cone and its pressure}

Write $x=r\theta$, $\theta\in\mathbb S^{N-1}$.
The change of variables absorbs the radial weight into the metric and measure.
Following \cite{CiraoloCorso2022}, introduce the cone metric and measure
\begin{equation}\label{eq-cone-data}
\begin{gathered}
 \rho=r^{h_a},\qquad \widehat N=\frac{N+a}{h_a}=p+\frac{N-p}{h_a}>N,\qquad
 \gamma=\frac{\widehat N-p}{p}=\frac{N-p}{p+a},\\
 g=\dd\rho^2+h_a^2\rho^2\dd\theta^2,\qquad
 \dd\mu=\rho^{\widehat N-1}\dd\rho\dd\theta,\qquad H=h_a^{-p}.
\end{gathered}
\end{equation}
Here $\dd\theta^2$ is the unit round metric on $\mathbb S^{N-1}$, $\dd\theta$ is its surface measure, and $o$ is the cone vertex. Write $w(\rho,\theta)=u(r\theta)$ and $D=\nabla_g$. For a vector field $Y$, weighted divergence is defined by
\[
\begin{gathered}
 \varphi=-(\widehat N-N)\log\rho,\qquad
 \dd\mu=h_a^{1-N}e^{-\varphi}\dd\operatorname{vol}_g,\\
 \operatorname{div}_\mu Y
 =\operatorname{div}_gY-\langle D\varphi,Y\rangle_g.
\end{gathered}
\]
The weak equation becomes
\begin{equation}\label{eq-cone-pde}
 -\Delta_{p,g,\mu}w:=-\operatorname{div}_\mu(|Dw|^{p-2}Dw)
 =H w^{m-1},\qquad m=\frac{p\widehat N}{\widehat N-p}.
\end{equation}
For any measurable region and its image under $\rho=r^{h_a}$,
\begin{equation}\label{eq-cone-energies}
 \int|\nabla u|^p\dd x=h_a^{p-1}\int|Dw|^p\dd\mu,
 \qquad \int|x|^au^m\dd x=h_a^{-1}\int w^m\dd\mu.
\end{equation}
These change-of-variable identities hold first on compact annuli, where all integrals are finite.

Define
\begin{equation}\label{eq-critical-pressure}
\begin{gathered}
 v=w^{-1/\gamma},\qquad X_v=|Dv|^{p-2}Dv,\qquad
 b=\frac{\widehat N(p-1)}p,\quad c=\gamma^{1-p}H,\\
 P:=\frac{b|Dv|^p+c}{v}=\operatorname{div}_\mu X_v.
\end{gathered}
\end{equation}
The divergence identity holds weakly on $\{\rho>0\}$: test \eqref{eq-cone-pde} with $v^b\phi$, $\phi\in C_c^\infty(\{\rho>0\})$. It holds classically on $\mathcal O:=\{\rho>0,\ |Dv|>0\}$; only there will we differentiate $P$.
Lemma~\ref{lem-critical-prior} implies, for $\rho\ge1$,
\begin{equation}\label{eq-cone-coarse}
 c_u\rho\le v\le C_u\rho^{p'},\qquad
 |Dv|\le C_u v/\rho,\qquad p'=\frac p{p-1}.
\end{equation}
At the vertex,
\begin{equation}\label{eq-pressure-apex}
 |Dv|\le C_u\left(r^{\alpha_1-h_a}+r^{h_a/(p-1)}\right)\longrightarrow0,
 \qquad P(o)=\frac c{v(o)}.
\end{equation}
Thus $v$ and $P$ extend continuously to $o$. At infinity, \eqref{eq-cone-coarse} gives
\[
 P\le C_u\frac{v^{p-1}}{\rho^p}+\frac cv\le C_u.
\]
Thus $P\in L^\infty$ and, since $v\ge c_u\rho\to\infty$, $v_0:=\min v>0$ is attained.

\begin{lemma}\label{lem-cone-pressure}
On $\mathcal O$, $D^2$ denotes the covariant Hessian, $DX_v(Y)=\nabla^g_YX_v$, and $I$ is the identity on the tangent space. Transposes denote $g$-adjoints, and matrix norms are the corresponding operator norms. Let
\[
 A=|Dv|^{p-2}\left(I+(p-2)\frac{Dv\otimes Dv}{|Dv|^2}\right).
\]
Then
\begin{equation}\label{eq-pressure-bochner}
 \operatorname{div}_\mu(v^{2-\widehat N}A DP)
 =\widehat N(p-1)v^{1-\widehat N}\mathcal D\ge0,
\end{equation}
where
\begin{equation}\label{eq-cone-defect}
 \mathcal D=\operatorname{tr}((DX_v)^2)
 +\operatorname{Ric}_{\varphi}(X_v,X_v)-\frac{P^2}{\widehat N}.
\end{equation}
Here $\operatorname{Ric}_g$ is the Ricci tensor of $g$, and $\operatorname{Ric}_{\varphi}=\operatorname{Ric}_g+D^2\varphi$.
\end{lemma}

\begin{proof}
With $g_T=g-\dd\rho\otimes\dd\rho$ the tangential metric, the effective Ricci tensor is
\begin{equation}\label{eq-cone-curvature}
\begin{aligned}
 \operatorname{Ric}_{\widehat N}
 &:=\operatorname{Ric}_{\varphi}-\frac{D\varphi\otimes D\varphi}{\widehat N-N}
 =\frac\kappa{\rho^2}g_T,\\
 \kappa&=\frac{N-2}{h_a^2}-(\widehat N-2),\qquad
 h_a^2\kappa=(1-h_a)\bigl(N-2+(p-2)h_a\bigr)>0.
\end{aligned}
\end{equation}
The last sign uses $N\ge3$.

On $\mathcal O$, $v$ is smooth and $DX_v=A D^2v$. Put $s=\operatorname{div}_gX_v$, $t=\langle D\varphi,X_v\rangle$, so that $P=s-t$. Since $A=A^T>0$, the Newton inequality \cite[Eq.~(3.11)]{CFR2020} gives
\[
 A^{-1/2}(DX_v)A^{1/2}=A^{1/2}D^2vA^{1/2},\qquad
 \operatorname{tr}((DX_v)^2)\ge\frac{s^2}{N}.
\]
Therefore
\begin{equation}\label{eq-defect-squares}
 \mathcal D=
 \left[\operatorname{tr}((DX_v)^2)-\frac{s^2}{N}\right]
 +\frac{((\widehat N-N)s+Nt)^2}{N\widehat N(\widehat N-N)}
 +\operatorname{Ric}_{\widehat N}(X_v,X_v)\ge0.
\end{equation}
Here $\widehat N>N$ and \eqref{eq-cone-curvature} make all three terms nonnegative.

Let $E=DX_v-(P/\widehat N)I$. Differentiating $vP=b|Dv|^p+c$ gives
\begin{equation}\label{eq-pressure-gradient}
 E^{T}Dv=\frac v{\widehat N} DP,\qquad
 EX_v=\frac{v}{\widehat N(p-1)}A DP.
\end{equation}
The covariant divergence formula gives
\[
 \operatorname{div}_\mu(DX_v\,X_v)
 =\langle X_v,DP\rangle+\operatorname{tr}((DX_v)^2)
  +\operatorname{Ric}_{\varphi}(X_v,X_v),
\]
and hence
\[
 \operatorname{div}_\mu(EX_v)=\frac{\widehat N-1}{\widehat N}\langle X_v,DP\rangle+\mathcal D.
\]
Since $A Dv=(p-1)X_v$, \eqref{eq-pressure-gradient} gives
\[
 \langle Dv,EX_v\rangle=\frac v{\widehat N}\langle X_v,DP\rangle.
\]
Hence
\begin{align*}
 \operatorname{div}_\mu(v^{1-\widehat N}EX_v)
 &=v^{1-\widehat N}\operatorname{div}_\mu(EX_v)
   +(1-\widehat N)v^{-\widehat N}\langle Dv,EX_v\rangle\\
 &=v^{1-\widehat N}\mathcal D
 =\frac1{\widehat N(p-1)}\operatorname{div}_\mu(v^{2-\widehat N}A DP),
\end{align*}
which is \eqref{eq-pressure-bochner}.
\end{proof}

\begin{lemma}\label{lem-sharp-critical-pressure}
With $v_0=\min v$ and $k_0=c/v_0$,
\begin{equation}\label{eq-sharp-critical-pressure}
 P\le k_0,\qquad |Dv|^p\le\frac{c}{b v_0}(v-v_0).
\end{equation}
\end{lemma}

\begin{proof}
For $k>k_0$, \eqref{eq-critical-pressure} gives
\[
 \left.P\right|_{\{Dv=0\}\cup\{o\}}=\frac cv\le k_0,\qquad
 \Omega_k:=\{P>k\},\qquad
 \overline{\Omega_k}\subset\mathcal O.
\]
Thus compactly supported tests in $\Omega_k$ require no extension across $\{Dv=0\}$.

For $R>2$, choose a radial cutoff $\eta_R$, equal to $1$ on $\{\rho\le R\}$ and $0$ on $\{\rho\ge R^2\}$, with
\[
 |D\eta_R|\le\frac{C}{\rho\log R}.
\]
For $\rho\ge1$ on $\Omega_k$, \eqref{eq-cone-coarse} gives, if $p\ge2$,
\[
 v^{2-\widehat N}\norm A\le C\rho^{2-p}v^{p-\widehat N}\le C\rho^{2-\widehat N}.
\]
If $1<p<2$, then
\[
 b|Dv|^p=vP-c\ge(k-k_0)v,
\]
so
\[
 v^{2-\widehat N}\norm A\le C_kv^{3-\widehat N-2/p}
 \le C_k\rho^{3-\widehat N-2/p}\le C_k\rho^{2-\widehat N}.
\]
Here $\widehat N>N\ge3$ makes $3-\widehat N-2/p<0$. Therefore
\begin{equation}\label{eq-pressure-log-capacity}
 \int_{\Omega_k}v^{2-\widehat N}\langle A D\eta_R,D\eta_R\rangle\dd\mu
 \le\frac{C_k}{(\log R)^2}\int_R^{R^2}\frac{\dd\rho}{\rho}
 =\frac{C_k}{\log R}\longrightarrow0.
\end{equation}
The inclusion $\overline{\Omega_k}\subset\mathcal O$ makes $\eta_R^2(P-k)_+$ compactly supported and Lipschitz in $\mathcal O$. Smooth approximation makes it admissible in \eqref{eq-pressure-bochner}. Testing and Cauchy--Schwarz yield
\begin{align*}
 &\int\eta_R^2v^{2-\widehat N}
 \langle A D(P-k)_+,D(P-k)_+\rangle\dd\mu\\
 &\qquad\le4\norm{(P-k)_+}_{L^\infty}^2
 \int_{\Omega_k}v^{2-\widehat N}\langle A D\eta_R,D\eta_R\rangle\dd\mu
 \longrightarrow0.
\end{align*}
Hence $DP=0$ on $\Omega_k$. If a connected component were nonempty, $P$ would be constant there with value greater than $k$. Continuity then makes that component closed as well as open in the connected punctured cone. This contradicts $P<k$ near the vertex. Thus $P\le k$; let $k\downarrow k_0$ and use \eqref{eq-critical-pressure}.
\end{proof}

\subsection{The center and the Picone comparison}

For the explicit profile $Z$ below, the remaining steps are
\[
 P\le\frac c{v_0}
 \quad\Longrightarrow\quad v(o)=v_0
 \quad\Longrightarrow\quad u\ge Z
 \quad\xRightarrow{\text{Picone}}\quad u=Z.
\]

\begin{proof}[Proof of Theorem~\ref{thm-negative-critical}]

Let $y_0$ be a minimum point of $v$, and set
\[
 B=\frac1{p'}\left(\frac{k_0}{\widehat N}\right)^{1/(p-1)},\qquad
 \widehat N(p'B)^{p-1}=k_0.
\]
For $\eps>0$, \eqref{eq-sharp-critical-pressure} gives
\[
 \left|D(v-v_0+\eps)^{1/p'}\right|
 \le\frac1{p'}\left(\frac c{b v_0}\right)^{1/p}
 =B^{1/p'}.
\]
Integrate along a unit-speed curve from $y_0$ and let $\eps\downarrow0$. At arclength $s$,
\begin{equation}\label{eq-critical-distance-growth}
 v-v_0\le Bs^{p'},\qquad |X_v|=|Dv|^{p-1}\le\frac{k_0}{\widehat N}s.
\end{equation}

Suppose $y_0\ne o$ and let $E_s=B_s^g(y_0)$ be a sufficiently small geodesic ball, with outward metric unit normal $\nu$ and weighted surface measure $\dd\sigma_\mu$. The weak divergence identity \eqref{eq-critical-pressure} gives, for almost every $s$,
\begin{equation}\label{eq-critical-local-flux}
 \frac1{\mu(E_s)}\int_{E_s}P\dd\mu
 =\frac1{\mu(E_s)}\int_{\partial E_s}\langle X_v,\nu\rangle\dd\sigma_\mu
 \le\frac{k_0}{\widehat N}\frac{s\,\sigma_\mu(\partial E_s)}{\mu(E_s)}.
\end{equation}
At this smooth point, $P(y_0)=k_0$ and the measure has a smooth positive density. Hence
\[
 \frac1{\mu(E_s)}\int_{E_s}P\dd\mu\longrightarrow k_0,
 \qquad
 \frac{s\,\sigma_\mu(\partial E_s)}{\mu(E_s)}\longrightarrow N.
\]
Thus \eqref{eq-critical-local-flux} would give
\[
 0<k_0\le\frac N{\widehat N} k_0<k_0.
\]
Hence $y_0=o$ and $v(o)=v_0$.

Apply \eqref{eq-critical-distance-growth} along radial curves from $o$. Since $\rho=r^{h_a}$, we obtain
\begin{equation}\label{eq-critical-comparison-function}
 u(x)\ge Z(x):=\left(v_0+B|x|^k\right)^{-\gamma},
 \qquad k=\frac{p+a}{p-1}.
\end{equation}
The normalization is
\begin{equation}\label{eq-critical-comparison-constant}
 \widehat N(\gamma Bp')^{p-1}v_0=\gamma^{p-1}k_0v_0=H.
\end{equation}
By \eqref{eq-bubble},
\[
 Z=\lambda^\beta U(\lambda x),\qquad
 \lambda=(B/v_0)^{1/k},\qquad
 -\Delta_pZ=|x|^aZ^{m-1}.
\]
The bounds
\[
 Z=O(r^{-(N-p)/(p-1)}),\qquad
 |\nabla Z|=O(r^{-(N-1)/(p-1)})\quad(r\to\infty)
\]
and the origin estimates following \eqref{eq-bubble} give
\begin{equation}\label{eq-critical-Z-energy}
 Z\in\D,\qquad
 \int|\nabla Z|^p\dd x=\int|x|^aZ^m\dd x<\infty.
\end{equation}
The identity follows by integrating the radial equation against $Z$; the flux terms vanish at $0$ and infinity.

Apply \eqref{eq-picone} with $V=|x|^au^{m-p}$ and $\phi=Z\in\D$. Since $u\ge Z$ and $m>p$,
\begin{equation}\label{eq-critical-Picone-chain}
 0\le\int|x|^a\bigl(u^{m-p}-Z^{m-p}\bigr)Z^p\dd x
 \le\int|\nabla Z|^p\dd x-\int|x|^aZ^m\dd x=0.
\end{equation}
Therefore
\[
 u=Z=\lambda^\beta U(\lambda\cdot)\in\D,\qquad
 \int|\nabla u|^p=\int|x|^au^m<\infty.
\]
\end{proof}

\begin{remark}
The restrictions enter through
\[
\begin{aligned}
 N=2&\quad\Longrightarrow\quad
 h_a^2\kappa=(1-h_a)(p-2)h_a<0,\\
 a=0&\quad\Longrightarrow\quad m=p^*,\\
 a>0&\quad\Longrightarrow\quad m>p^*.
\end{aligned}
\]
Thus the curvature argument requires $N\ge3$, and the subcritical blow-up contradiction requires $a<0$. For $a>0$, known nonradial branches \cite{DDGL2026} also exclude classification by radial profiles alone.
\end{remark}


\begin{thebibliography}{99}
\setlength{\itemsep}{0pt}
\raggedright

\bibitem{AH1998}
W. Allegretto and Y. X. Huang,
A Picone's identity for the $p$-Laplacian and applications,
Nonlinear Anal. \textbf{32} (1998), no. 7, 819--830.
\doi{10.1016/S0362-546X(97)00530-0}.


\bibitem{BV1989}
M.-F. Bidaut-V\'eron,
Local and global behavior of solutions of quasilinear equations of Emden--Fowler type,
Arch. Rational Mech. Anal. \textbf{107} (1989), no. 4, 293--324.
\doi{10.1007/BF00251552}.

\bibitem{BGH2001}
M.-F. Bidaut-V\'eron and M. Garc\'ia-Huidobro,
Regular and singular solutions of a quasilinear equation with weights,
Asymptot. Anal. \textbf{28} (2001), no. 2, 115--150.

\bibitem{BidautPohozaev2001}
M.-F. Bidaut-V\'eron and S. I. Pohozaev,
Nonexistence results and estimates for some nonlinear elliptic problems,
J. Anal. Math. \textbf{84} (2001), 1--49.
\doi{10.1007/BF02788105}.

\bibitem{Brock1995}
F. Brock,
Continuous Steiner-symmetrization,
Math. Nachr. \textbf{172} (1995), 25--48.
\doi{10.1002/mana.19951720104}.

\bibitem{Brock2000}
F. Brock,
Continuous rearrangement and symmetry of solutions of elliptic problems,
Proc. Indian Acad. Sci. Math. Sci. \textbf{110} (2000), no. 2, 157--204.
\doi{10.1007/BF02829490}.

\bibitem{Brock}
F. Brock,
Symmetry for a general class of overdetermined elliptic problems,
Nonlinear Differ. Equ. Appl. NoDEA \textbf{23} (2016), article 36.
\doi{10.1007/s00030-016-0390-1}.

\bibitem{CiraoloCorso2022}
G. Ciraolo and R. Corso,
Symmetry for positive critical points of Caffarelli--Kohn--Nirenberg inequalities,
Nonlinear Anal. \textbf{216} (2022), article 112683.
\doi{10.1016/j.na.2021.112683}.

\bibitem{CFR2020}
G. Ciraolo, A. Figalli, and A. Roncoroni,
Symmetry results for critical anisotropic $p$-Laplacian equations in convex cones,
Geom. Funct. Anal. \textbf{30} (2020), no. 3, 770--803.
\doi{10.1007/s00039-020-00535-3}.

\bibitem{DDGL2026}
W. Dai, L. Duan, C. Gui, and Y. Li,
Nonradial solutions for the critical quasi-linear H\'enon equation involving $p$-Laplacian in $\R^N$,
Proc. Lond. Math. Soc. \textbf{132} (2026), e70148.
\doi{10.1112/plms.70148}.

\bibitem{DMS}
M. Degiovanni, A. Musesti, and M. Squassina,
On the regularity of solutions in the Pucci--Serrin identity,
Calc. Var. Partial Differential Equations \textbf{18} (2003), no. 3, 317--334.
\doi{10.1007/s00526-003-0208-y}.

\bibitem{DengZhang2026}
S. Deng and B. Zhang,
Liouville theorems and sharp homogeneity defects for weighted potential systems,
preprint (2026).
\href{https://arxiv.org/abs/2610.31167v1}{arXiv:2610.31167v1}.

\bibitem{DiBenedetto1983}
E. DiBenedetto,
$C^{1+\alpha}$ local regularity of weak solutions of degenerate elliptic equations,
Nonlinear Anal. \textbf{7} (1983), no. 8, 827--850.
\doi{10.1016/0362-546X(83)90061-5}.

\bibitem{DZ}
G. Du and S. Zhou,
Local and global properties of $p$-Laplace H\'enon equation,
J. Math. Anal. Appl. \textbf{535} (2024), no. 2, article 128148.
\doi{10.1016/j.jmaa.2024.128148}.

\bibitem{GarciaMelian2017}
J. Garc\'ia-Meli\'an,
Nonexistence of positive solutions for H\'enon equation,
preprint (2017).
\href{https://arxiv.org/abs/1703.04353}{arXiv:1703.04353}.

\bibitem{GhoussoubYuan2000}
N. Ghoussoub and C. Yuan,
Multiple solutions for quasi-linear PDEs involving the critical Sobolev and Hardy exponents,
Trans. Amer. Math. Soc. \textbf{352} (2000), no. 12, 5703--5743.
\doi{10.1090/S0002-9947-00-02560-5}.

\bibitem{GS1981}
B. Gidas and J. Spruck,
Global and local behavior of positive solutions of nonlinear elliptic equations,
Comm. Pure Appl. Math. \textbf{34} (1981), no. 4, 525--598.
\doi{10.1002/cpa.3160340406}.

\bibitem{GigaNgo2022}
Y. Giga and Q. A. Ng\^o,
Exhaustive existence and non-existence results for Hardy--H\'enon equations in $\R^n$,
Partial Differ. Equ. Appl. \textbf{3} (2022), article 81.
\doi{10.1007/s42985-022-00190-3}.


\bibitem{Henon1973}
M. H\'enon,
Numerical experiments on the stability of spherical stellar systems,
Astronom. Astrophys. \textbf{24} (1973), 229--238.
\href{https://ui.adsabs.harvard.edu/abs/1973A\%26A....24..229H/abstract}{NASA ADS}.

\bibitem{KNY1990}
N. Kawano, W.-M. Ni, and S. Yotsutani,
A generalized Pohozaev identity and its applications,
J. Math. Soc. Japan \textbf{42} (1990), no. 3, 541--564.
\doi{10.2969/jmsj/04230541}.

\bibitem{KM1994}
T. Kilpel\"ainen and J. Mal\'y,
The Wiener test and potential estimates for quasilinear elliptic equations,
Acta Math. \textbf{172} (1994), no. 1, 137--161.
\doi{10.1007/BF02392793}.

\bibitem{Kilpelainen}
T. Kilpel\"ainen,
$p$-Laplacian type equations involving measures,
in \emph{Proceedings of the International Congress of Mathematicians},
Vol. III, Beijing, 2002, Higher Education Press, Beijing, 2002, 167--176.

\bibitem{KZ2002}
T. Kilpel\"ainen and X. Zhong,
Removable sets for continuous solutions of quasilinear elliptic equations,
Proc. Amer. Math. Soc. \textbf{130} (2002), no. 6, 1681--1688.
\doi{10.1090/S0002-9939-01-06237-2}.

\bibitem{LeLe2023}
P. Le and D. H. T. Le,
Classification of positive solutions to $p$-Laplace equations with critical Hardy--Sobolev exponent,
Nonlinear Anal. Real World Appl. \textbf{74} (2023), article 103949.
\doi{10.1016/j.nonrwa.2023.103949}.
		
\bibitem{Le2026Anisotropic}
P. Le,
Critical anisotropic $p$-Laplace equations without finite energy:
classification and a Schoen-type Harnack inequality,
preprint (2026).
\href{https://arxiv.org/abs/2610.00337}{arXiv:2610.00337}.

\bibitem{LiLuZhang2026Fractional}
K. Li, H. Lu, and Z. Zhang,
The fractional H\'enon--Lane--Emden conjecture and some critical nonexistence phenomena,
preprint (2026).
\href{https://arxiv.org/abs/2610.22802}{arXiv:2610.22802}.

\bibitem{LiZhang2019Monotonicity}
K. Li and Z. Zhang,
A monotonicity theorem and its applications to weighted elliptic equations,
Sci. China Math. \textbf{62} (2019), 1925--1934.
\doi{10.1007/s11425-018-9414-8}.

\bibitem{Lindqvist2019}
P. Lindqvist,
\emph{Notes on the Stationary $p$-Laplace Equation},
SpringerBriefs in Mathematics, Springer, Cham, 2019.
\doi{10.1007/978-3-030-14501-9}.

\bibitem{LuZhang2026Disk}
H. Lu and Z. Zhang,
New result on Lane--Emden conjecture with exponents in a disk,
preprint (2026).
\href{https://arxiv.org/abs/2610.09526}{arXiv:2610.09526}.

\bibitem{MW2026}
L. Mei and J. Wei,
A proof of the fractional Lane--Emden conjecture,
preprint (2026).
\href{https://arxiv.org/abs/2610.10188}{arXiv:2610.10188}.

\bibitem{MP2001}
E. Mitidieri and S. I. Pohozaev,
A priori estimates and the absence of solutions of nonlinear partial differential equations and inequalities,
Proc. Steklov Inst. Math. \textbf{234} (2001), no. 3, 1--362.

\bibitem{Ni1982}
W.-M. Ni,
A nonlinear Dirichlet problem on the unit ball and its applications,
Indiana Univ. Math. J. \textbf{31} (1982), no. 6, 801--807.
\doi{10.1512/iumj.1982.31.31056}.

\bibitem{OpenAI2026}
OpenAI,
The subcritical H\'enon--Lane--Emden conjecture,
research preprint, 24 September 2026.
\href{https://github.com/openai/math/tree/main/preprints/The-Subcritical-Henon-Lane-Emden-Conjecture-September-24-2026}{Public manuscript in the OpenAI math repository}.

\bibitem{PanShang2026}
Y. Pan and Y. Shang,
Classification of positive solution to the $p$-Laplace equation with critical Sobolev--Hardy exponent,
Complex Var. Elliptic Equ. (2026), published online 4 May 2026.
\doi{10.1080/17476933.2026.2640582}.

\bibitem{PD}
Q. H. Phan and A. T. Duong,
Liouville-type theorems for a quasilinear elliptic equation of the H\'enon-type,
Nonlinear Differ. Equ. Appl. NoDEA \textbf{22} (2015), no. 6, 1817--1829.
\doi{10.1007/s00030-015-0345-y}.

\bibitem{PhanSouplet2012}
Q. H. Phan and P. Souplet,
Liouville-type theorems and bounds of solutions of Hardy--H\'enon equations,
J. Differential Equations \textbf{252} (2012), no. 3, 2544--2562.
\doi{10.1016/j.jde.2011.09.022}.

\bibitem{PQS2007}
P. Pol\'a\v cik, P. Quittner, and P. Souplet,
Singularity and decay estimates in superlinear problems via Liouville-type theorems. I. Elliptic equations and systems,
Duke Math. J. \textbf{139} (2007), no. 3, 555--579.
\doi{10.1215/S0012-7094-07-13935-8}.

\bibitem{PS}
P. Pucci and R. Servadei,
Existence, non-existence and regularity of radial ground states for $p$-Laplacian equations with singular weights,
Ann. Inst. H. Poincar\'e Anal. Non Lin\'eaire \textbf{25} (2008), no. 3, 505--537.
\doi{10.1016/j.anihpc.2007.02.004}.

\bibitem{PS2010}
P. Pucci and R. Servadei,
Nonexistence for $p$-Laplace equations with singular weights,
Commun. Pure Appl. Anal. \textbf{9} (2010), no. 5, 1421--1438.
\doi{10.3934/cpaa.2010.9.1421}.

\bibitem{Serrin64}
J. Serrin,
Local behavior of solutions of quasi-linear equations,
Acta Math. \textbf{111} (1964), 247--302.
\doi{10.1007/BF02391014}.

\bibitem{SZ}
J. Serrin and H. Zou,
Cauchy--Liouville and universal boundedness theorems for quasilinear elliptic equations and inequalities,
Acta Math. \textbf{189} (2002), no. 1, 79--142.
\doi{10.1007/BF02392645}.

\bibitem{Tolksdorf}
P. Tolksdorf,
Regularity for a more general class of quasilinear elliptic equations,
J. Differential Equations \textbf{51} (1984), no. 1, 126--150.
\doi{10.1016/0022-0396(84)90105-0}.

\bibitem{Zhang2026Critical}
Yi Ru-Ya Zhang,
Critical $p$-Laplace equations with monotone coefficients:
Liouville classification and a Schoen-type Harnack inequality,
preprint (2026).
\href{https://arxiv.org/abs/2608.09113}{arXiv:2608.09113}.

\bibitem{zzt}
Z. Zhang,
\emph{Variational, Topological, and Partial Order Methods with Their Applications},
Developments in Mathematics, vol.~29, Springer, Heidelberg, 2013.
\doi{10.1007/978-3-642-30709-6}.






\end{thebibliography}
\end{document}